\documentclass[11pt,a4paper,thmsb,ukenglish,dvips]{article}
\usepackage{amsfonts}
\usepackage{amssymb}
\usepackage{amsmath}
\usepackage{chicago}
\usepackage{epsf}
\usepackage{theorem}

\newtheorem{theorem}{Theorem}

\newtheorem{corollary}{Corollary}

\newtheorem{example}{Example}

\newtheorem{lemma}{Lemma}

\newtheorem{remark}{Remark}

\newenvironment{assumption}[1][Assumption]{\noindent\textbf{#1} \it}{\
}
\input{tcilatex}

\begin{document}

\title{Limit theorems for bipower variation \\
in financial econometrics}
\author{\textsc{Ole E. Barndorff-Nielsen } \\
%EndAName
\textit{Department of Mathematical Sciences,}\\
\textit{University of Aarhus, Ny Munkegade, DK-8000 Aarhus C, Denmark} \\
\texttt{oebn@imf.au.dk} \and \textsc{Svend Erik Graversen} \\
%EndAName
\textit{Department of Mathematical Sciences,}\\
\textit{University of Aarhus, Ny Munkegade, DK-8000 Aarhus C, Denmark}\\
\texttt{matseg@imf.au.dk} \and \textsc{Jean Jacod} \\
%EndAName
\textit{Laboratoire de Probabilit\'{e}s et Mod\`{e}les Al\'{e}atoires (CNRS
UMR 7599)}\\
\textit{Universit\'{e} Pierre et Marie Curie,}\\
\textit{4 Place Jussieu, 75252 Paris Cedex 05, France } \\
\texttt{jj@ccr.jussieu.fr} \and \textsc{Neil Shephard} \\
%EndAName
\textit{Nuffield College, University of Oxford, Oxford OX1 1NF, UK}\\
\texttt{neil.shephard@nuf.ox.ac.uk}}
\date{}
\maketitle

\begin{abstract}
In this paper we provide an asymptotic analysis of generalised bipower
measures of the variation of price processes in financial economics. These
measures encompass the usual quadratic variation, power variation and
bipower variations which have been highlighted in recent years in financial
econometrics. The analysis is carried out under some rather general Brownian
semimartingale assumptions, which allow for standard leverage effects.
\end{abstract}
\bigskip
\bigskip


\medskip \baselineskip=20pt

\newpage

\section{Introduction}

In this paper we discuss the limiting theory for a novel, unifying class of
non-parametric measures of the variation of financial prices. The theory
covers commonly used estimators of variation such as realised volatility,
but it also encompasses more recently suggested quantities like realised
power variation and realised bipower variation. We considerably strengthen
existing results on the latter two quantities, deepening our understanding
and unifying their treatment. We will outline the proofs of these theorems,
referring for the very technical, detailed formal proofs of the general
results to a companion probability theory paper \cite%
{BarndorffNielsenGraversenJacodPodolskyShephard(04shiryaev)}. Our emphasis
is on exposition, explaining where the results come from and how they sit
within the econometrics literature.

Our theoretical development is motivated by the advent of complete records
of quotes or transaction prices for many financial assets. Although market
microstructure effects (e.g. discreteness of prices, bid/ask bounce,
irregular trading etc.) mean that there is a mismatch between asset pricing
theory based on semimartingales and the data at very fine time intervals it
does suggest the desirability of establishing an asymptotic distribution
theory for estimators as we use more and more highly frequent observations.
Papers which directly model the impact of market microstructure noise on
realised variance include \cite{BandiRussell(03)}, \cite{HansenLunde(03)},
\cite{ZhangMyklandAitSahalia(03)}, \cite%
{BarndorffNielsenHansenLundeShephard(04)} and \cite{Zhang(04)}. Related work
in the probability literature on the impact of noise on discretely observed
diffusions can be found in \cite{GloterJacod(01a)} and \cite%
{GloterJacod(01b)}, while \cite{DelattreJacod(97)} report results on the
impact of rounding on sums of functions of discretely observed diffusions.
In this paper we ignore these effects.

Let the $d$-dimensional vector of the log-prices of a set of assets follow
the process
\begin{equation*}
Y=\left( Y^{1},...,Y^{d}\right) ^{\prime }.
\end{equation*}
At time $t\geq 0$ we denote the log-prices as $Y_{t}$. Our aim is to
calculate measures of the variation of the price process (e.g. realised
volatility) over discrete time intervals (e.g. a day or a month). Without
loss of generality we can study the mathematics of this by simply looking at
what happens when we have $n$ high frequency observations on the time
interval $t=0$ to $t=1$ and study what happens to our measures of variation
as $n\rightarrow \infty $ (e.g., for introductions to this, \cite%
{BarndorffNielsenShephard(02realised)}). In this case returns will be
measured over intervals of length $n^{-1}$ as
\begin{equation}
\Delta _{i}^{n}Y=Y_{i/n}-Y_{(i-1)/n},\quad i=1,2,...,n,  \label{return}
\end{equation}%
where $n$ is a positive integer.

We will study the behaviour of the realised generalised bipower variation
process
\begin{equation}
\frac{1}{n}\sum_{i=1}^{\left\lfloor nt\right\rfloor }g(\sqrt{n}~\Delta
_{i}^{n}Y)h(\sqrt{n}~\Delta _{i+1}^{n}Y),  \label{RGBP}
\end{equation}%
as $n$ becomes large and where $g$ and $h$ are two given, matrix functions
of dimensions $d_{1}\times d_{2}$ and $d_{2}\times d_{3}$ respectively,
whose elements have at most polynomial growth. Here $\left\lfloor
x\right\rfloor $ denotes the largest integer less than or equal to $x$.

Although (\ref{RGBP}) looks initially rather odd, in fact most of the
non-parametric volatility measures used in financial econometrics fall
within this class (a measure not included in this setup is the range
statistic studied in, for example, \cite{Parkinson(80)}). Here we give an
extensive list of examples and link them to the existing literature. More
detailed discussion of the literature on the properties of these special
cases will be given later.

\begin{example}
\label{Example: 1}\textbf{(a)} Suppose $g(y)=\left( y^{j}\right) ^{2}$ and $%
h(y)=1$, then (\ref{RGBP}) becomes%
\begin{equation*}
\sum_{i=1}^{\left\lfloor nt\right\rfloor }\left( \Delta _{i}^{n}Y^{j}\right)
^{2},
\end{equation*}%
which is called the realised quadratic variation process of $Y^{j}$ in
econometrics, e.g. \cite{Jacod(94)}, \cite{JacodProtter(98)}, \cite%
{BarndorffNielsenShephard(02realised)}, \cite%
{BarndorffNielsenShephard(04multi)} and \cite{MyklandZhang(05)}. The
increments of this quantity, typically calculated over a day or a week, are
often called the realised variances in financial economics and have been
highlighted by \cite{AndersenBollerslevDieboldLabys(01)} and \cite%
{AndersenBollerslevDiebold(05)} in the context of volatility measurement and
forecasting.

\noindent \textbf{(b)} Suppose $g(y)=yy^{\prime }$ and $h(y)=I$, then (\ref%
{RGBP}) becomes, after some simplification,
\begin{equation*}
\sum_{i=1}^{\left\lfloor nt\right\rfloor }\left( \Delta _{i}^{n}Y\right)
\left( \Delta _{i}^{n}Y\right) ^{\prime }.
\end{equation*}%
\newline
This is the realised covariation process. It has been studied by \cite%
{JacodProtter(98)}, \cite{BarndorffNielsenShephard(04multi)} and \cite%
{MyklandZhang(05)}. \cite{AndersenBollerslevDieboldLabys(03model)} study the
increments of this process to produce forecast distributions for vectors of
returns. \

\noindent \textbf{(c)} Suppose $g(y)=\left\vert y^{j}\right\vert ^{r}$ for $%
r>0$ and $h(y)=1$, then (\ref{RGBP}) becomes
\begin{equation*}
n^{-1+r/2}\sum_{i=1}^{\left\lfloor nt\right\rfloor }\left\vert \Delta
_{i}^{n}Y^{j}\right\vert ^{r},
\end{equation*}%
which is called the realised $r$-th order power variation. When $r$ is an
integer it has been studied from a probabilistic viewpoint by \cite%
{Jacod(94)} while \cite{BarndorffNielsenShephard(03bernoulli)} look at the
econometrics of the case where $r>0$. The increments of these types of high
frequency volatility measures have been informally used in the financial
econometrics literature for some time when $r=1$, but until recently without
a strong understanding of their properties. Examples of their use include
\cite{Schwert(90JB)}, \cite{AndersenBollerslev(98)} and \cite%
{AndersenBollerslev(97jef)}, while they have also been informally discussed
by \cite[pp. 349--350]{Shiryaev(99)}\ and \cite{MaheswaranSims(93)}.
Following the work by \cite{BarndorffNielsenShephard(03bernoulli)}, \cite%
{GhyselsSantaClaraValkoanov(04)} and \cite{ForsbergGhysels(04)} have
successfully used realised power variation as an input into volatility
forecasting competitions.

\noindent \textbf{(d)} Suppose $g(y)=\left\vert y^{j}\right\vert ^{r}$ and $%
h(y)=\left\vert y^{j}\right\vert ^{s}$ for $r,s>0$, then (\ref{RGBP})
becomes
\begin{equation*}
n^{-1+(r+s)/2}\sum_{i=1}^{\left\lfloor nt\right\rfloor }\left\vert \Delta
_{i}^{n}Y^{j}\right\vert ^{r}\left\vert \Delta _{i+1}^{n}Y^{j}\right\vert
^{s},
\end{equation*}%
which is called the realised $r,s$-th order bipower variation process. This
measure of variation was introduced by \cite{BarndorffNielsenShephard(04jfe)}%
, while a more formal discussion of its behaviour in the $r=s=1$ case was
developed by \cite{BarndorffNielsenShephard(03test)}. These authors'
interest in this quantity was motivated by its virtue of being resistant to
finite activity jumps so long as $\max (r,s)<2$. Recently \cite%
{BarndorffNielsenShephardWinkel(04)} and \cite{Woerner(04power)} have
studied how these results on jumps extend to infinite activity processes,
while \cite{CorradiDistaso(04)} have used these statistics to test the
specification of parametric volatility models.

\noindent \textbf{(e)} Suppose
\begin{equation*}
g(y)=\left(
\begin{array}{cc}
\left\vert y^{j}\right\vert & 0 \\
0 & \left( y^{j}\right) ^{2}%
\end{array}%
\right) ,\quad h(y)=\left(
\begin{array}{c}
\left\vert y^{j}\right\vert \\
1%
\end{array}%
\right) .
\end{equation*}%
Then (\ref{RGBP}) becomes,%
\begin{equation*}
\left(
\begin{array}{c}
\displaystyle\sum_{i=1}^{\left\lfloor nt\right\rfloor }\left\vert \Delta
_{i}^{n}Y^{j}\right\vert \left\vert \Delta _{i+1}^{n}Y^{j}\right\vert \\
\displaystyle\sum_{i=1}^{\left\lfloor nt\right\rfloor }\left( \Delta
_{i}^{n}Y^{j}\right) ^{2}%
\end{array}%
\right) .
\end{equation*}%
\cite{BarndorffNielsenShephard(03test)} used the joint behaviour of the
increments of these two statistics to test for jumps in price processes. \
\cite{HuangTauchen(03)} have empirically studied the finite sample
properties of these types of jump tests. \cite%
{AndersenBollerslevDiebold(03bipower)} \ and \cite{ForsbergGhysels(04)} use
bipower variation as an input into volatility forecasting. \
\end{example}

We will derive the probability limit of (\ref{RGBP}) under a general
Brownian semimartingale, the workhorse process of modern continuous time
asset pricing theory. Only the case of realised quadratic variation, where
the limit is the usual quadratic variation QV (defined for general
semimartingales), has been previously been studied under such wide
conditions. Further, under some stronger but realistic conditions, we will
derive a limiting distribution theory for (\ref{RGBP}), so extending a
number of results previously given in the literature on special cases of
this framework.

The outline of this paper is as follows. Section 2 contains a detailed
listing of the assumptions used in our analysis. Section 3 gives a statement
of a weak law of large numbers for these statistics and the corresponding
central limit theory is presented in Section 4. Extensions of the results to
higher order variations is briefly indicated in Section 5. Section 6
illustrates the theory by discussing how it gives rise to tests for jumps in
the price processes, using bipower and tripower variation. The corresponding
literature which discusses various special cases of these results is also
given in these sections. Section 8 concludes, while there is an Appendix
which provides an outline of the proofs of the results discussed in this
paper. For detailed, quite lengthy and highly technical formal proofs we
refer to our companion probability theory paper \cite%
{BarndorffNielsenGraversenJacodPodolskyShephard(04shiryaev)}.

\section{Notation and models}

We start with $Y$ on some filtered probability space $\left( \Omega ,%
\mathcal{F},\left( \mathcal{F}_{t}\right) _{t\geq 0},P\right) $. In most of
our analysis we will assume that $Y$ follows a $d$-dimensional Brownian
semimartingale (written $Y\in \mathcal{BSM}$). It is given in the following
statement.

\noindent \textbf{Assumption (H): }We have
\begin{equation}
Y_{t}=Y_{0}+\int_{0}^{t}a_{u}\mathrm{d}u+\int_{0}^{t}\sigma _{u-}\mathrm{d}%
W_{u},  \label{H}
\end{equation}%
where $W$ is a $d^{\prime }$-dimensional standard Brownian motion (BM), $a$
is a $d$-dimensional process whose elements are predictable and has locally
bounded sample paths, and the spot covolatility $d,d^{\prime }$-dimensional
matrix $\sigma $ has elements which have c\`{a}dl\`{a}g sample paths.

Throughout we will write
\begin{equation*}
\Sigma _{t}=\sigma _{t}\sigma _{t}^{\prime },
\end{equation*}%
the spot covariance matrix. Typically $\Sigma _{t}$ will be full rank, but
we do not assume that here. We will write $\Sigma _{t}^{jk}$ to denote the $%
j,k$-th element of $\Sigma _{t}$, while we write%
\begin{equation*}
\sigma _{j,t}^{2}=\Sigma _{t}^{jj}.
\end{equation*}

\begin{remark}
Due to the fact that $t\mapsto \sigma _{t}^{jk}$ is c\`{a}dl\`{a}g all
powers of $\sigma _{t}^{jk}$ are locally integrable with respect to the
Lebesgue measure. \ In particular then $\int_{0}^{t}\Sigma _{u}^{jj}\mathrm{d%
}u<\infty $ for all $t$ and $j$.
\end{remark}

\begin{remark}
Both $a$ and $\sigma $ can have, for example, jumps, intraday seasonality
and long-memory.
\end{remark}

\begin{remark}
The stochastic volatility (e.g. \cite{GhyselsHarveyRenault(96)} and \cite%
{Shephard(05)}) component of $Y$,
\begin{equation*}
\int_{0}^{t}\sigma _{u-}\mathrm{d}W_{u},
\end{equation*}%
is always a vector of local martingales each with continuous sample paths,
as $\int_{0}^{t}\Sigma _{u}^{jj}\mathrm{d}u<\infty $ for all $t$ and $j$.
All continuous local martingales with absolutely continuous quadratic
variation can be written in the form of a stochastic volatility process.
This result, which is due to \cite{Doob(53)}, is discussed in, for example,
\cite[p. 170--172]{KaratzasShreve(91)}. Using the Dambis-Dubins-Schwartz
Theorem, we know that the difference between the entire continuous local
martingale class and the SV class are the local martingales which have only
continuous, not absolutely continuous\footnote{%
An example of a continuous local martingale which has no SV representation
is a time-change Brownian motion where the time-change takes the form of the
so-called \textquotedblleft devil's staircase,\textquotedblright\ which is
continuous and non-decreasing but not absolutely continuous (see, for
example, \cite[Section 27]{Munroe(53)}). This relates to the work of, for
example, \cite{CalvetFisher(02)} on multifractals.}, QV. The drift $%
\int_{0}^{t}a_{u}\mathrm{d}u$ has elements which are absolutely continuous.
This assumption looks ad hoc, however if we impose a lack of arbitrage
opportunities and model the local martingale component as a SV process then
this property must hold (\cite[p. 3]{KaratzasShreve(98)} and \cite[p. 583]%
{AndersenBollerslevDieboldLabys(03model)}). Hence (\ref{H}) is a rather
canonical model in the finance theory of continuous sample path processes.
\end{remark}

We are interested in the asymptotic behaviour, for $n\rightarrow \infty $,
of the following volatility measuring process:
\begin{equation}
Y^{n}(g,h)_{t}=\frac{1}{n}\sum_{i=1}^{\left\lfloor nt\right\rfloor }g(\sqrt{n%
}~\Delta _{i}^{n}Y)h(\sqrt{n}~\Delta _{i+1}^{n}Y),  \label{XP}
\end{equation}%
where $g$ and $h$ are two given conformable matrix functions and recalling
the definition of $\Delta _{i}^{n}Y$ given in (\ref{return}).

\section{Law of large numbers}

To build a weak law of large numbers for $Y^{n}(g,h)_{t}$ we need to make
the pair $(g,h)$ satisfy the following assumption.

\noindent \textbf{Assumption (K):} All the elements of $f$ on $\mathbf{R}%
^{d} $ are continuous with at most polynomial growth.

This amounts to there being suitable constants $C>0$ and $p\geq 2$ such that
\begin{equation}
x\in \mathbf{R}^{d}\quad \Rightarrow \quad \left\Vert f(x)\right\Vert \leq
C(1+\Vert x\Vert ^{p}).  \label{G1}
\end{equation}

We also need the following notation.
\begin{equation*}
\rho _{\sigma }(g)=\mathrm{E}\left\{ g(X)\right\} ,\quad \text{where\quad }%
X|\sigma \sim N(0,\sigma \sigma ^{\prime }),
\end{equation*}%
and
\begin{equation*}
\rho _{\sigma }(gh)=\mathrm{E}\left\{ g(X)h(X)\right\} .
\end{equation*}

\begin{example}
\label{Example: second}\textbf{(a)} Let $g(y)=yy^{\prime }$ and $h(y)=I$,
then $\rho _{\sigma }(g)=\Sigma $ and $\rho _{\sigma }(h)=I$.

\noindent \textbf{(b)} Suppose $g(y)=\left\vert y^{j}\right\vert ^{r}$ then $%
\rho _{\sigma }(g)=\mu _{r}\sigma _{j}^{r}$, where $\sigma _{j}^{2}$ is the $%
j,j$-th element of $\Sigma $, $\mu _{r}=\mathrm{E}(\left\vert u\right\vert
^{r})$ and $u\sim N(0,1)$.
\end{example}

This setup is sufficient for the proof of Theorem 1.2 of \cite%
{BarndorffNielsenGraversenJacodPodolskyShephard(04shiryaev)}, which is
restated here.

\begin{theorem}
\label{TT1} Under (H) and assuming $g$ and $h$ satisfy (K) we have that
\begin{equation}
Y^{n}(g,h)_{t}~\rightarrow ~Y(g,h)_{t}:=\int_{0}^{t}\rho _{\sigma
_{u}}(g)\rho _{\sigma _{u}}(h)\mathrm{d}u,  \label{WLLN}
\end{equation}%
where the convergence is in probability, locally uniform in time. \newline
\end{theorem}

The result is quite clean as it is requires no additional assumptions on $Y$
and so is very close to dealing with the whole class of financially coherent
continuous sample path processes.

This Theorem covers a number of existing setups which are currently
receiving a great deal of attention as measures of variation in financial
econometrics. Here we briefly discuss some of the work which has studied the
limiting behaviour of these objects.

\begin{example}
\textbf{(Example \ref{Example: 1}(a) continued)}. Then $g(y)=\left(
y^{j}\right) ^{2}$ and $h(y)=1$, so (\ref{WLLN}) becomes%
\begin{equation*}
\sum_{i=1}^{\left\lfloor nt\right\rfloor }\left( \Delta _{i}^{n}Y^{j}\right)
^{2}\rightarrow ~\int_{0}^{t}\sigma _{j,u}^{2}\mathrm{d}u=[Y^{j}]_{t},
\end{equation*}%
the quadratic variation (QV) of $Y^{j}$. This well known result in
probability theory is behind much of the modern work on realised volatility,
which is compactly reviewed in \cite{AndersenBollerslevDiebold(05)}.

\noindent (\textbf{Example \ref{Example: 1}(b) continued}). As $%
g(y)=yy^{\prime }$ and $h(y)=I$, then
\begin{equation*}
\sum_{i=1}^{\left\lfloor nt\right\rfloor }\left( \Delta _{i}^{n}Y\right)
\left( \Delta _{i}^{n}Y\right) ^{\prime }\rightarrow ~\int_{0}^{t}\Sigma _{u}%
\mathrm{d}u=[Y]_{t},
\end{equation*}%
the well known multivariate version of QV.

\noindent \textbf{(Example \ref{Example: 1}(c) continued).} Then $%
g(y)=\left\vert y^{j}\right\vert ^{r}$ and $h(y)=1$ so
\begin{equation*}
n^{-1+r/2}\sum_{i=1}^{\left\lfloor nt\right\rfloor }\left\vert \Delta
_{i}^{n}Y^{j}\right\vert ^{r}\rightarrow ~\mu _{r}\int_{0}^{t}\sigma
_{j,u}^{r}\mathrm{d}u.
\end{equation*}%
This result is due to \cite{Jacod(94)} and \cite%
{BarndorffNielsenShephard(03bernoulli)}.

\noindent \textbf{(Example \ref{Example: 1}(d) continued).} Then $%
g(y)=\left\vert y^{j}\right\vert ^{r}$ and $h(y)=\left\vert y^{j}\right\vert
^{s}$ for $r,s>0$, so
\begin{equation*}
n^{-1+(r+s)/2}\sum_{i=1}^{\left\lfloor nt\right\rfloor }\left\vert \Delta
_{i}^{n}Y^{j}\right\vert ^{r}\left\vert \Delta _{i+1}^{n}Y^{j}\right\vert
^{s}\rightarrow ~\mu _{r}\mu _{s}\int_{0}^{t}\sigma _{j,u}^{r+s}\mathrm{d}u,
\end{equation*}%
a result due to \cite{BarndorffNielsenShephard(04jfe)}, who derived it under
stronger conditions than those used here.

\noindent \textbf{(Example \ref{Example: 1}(e) continued).} Then
\begin{equation*}
g(y)=\left(
\begin{array}{cc}
\left\vert y^{j}\right\vert & 0 \\
0 & \left( y^{j}\right) ^{2}%
\end{array}%
\right) ,\quad h(y)=\left(
\begin{array}{c}
\left\vert y^{j}\right\vert \\
1%
\end{array}%
\right) ,
\end{equation*}%
so
\begin{equation*}
\left(
\begin{array}{c}
\displaystyle\sum_{i=1}^{\left\lfloor nt\right\rfloor }\left\vert \Delta
_{i}^{n}Y^{j}\right\vert \left\vert \Delta _{i+1}^{n}Y^{j}\right\vert \\
\displaystyle\sum_{i=1}^{\left\lfloor nt\right\rfloor }\left( \Delta
_{i}^{n}Y^{j}\right) ^{2}%
\end{array}%
\right) \rightarrow \left( ~%
\begin{array}{c}
\mu _{1}^{2} \\
1%
\end{array}%
\right) \int_{0}^{t}\sigma _{j,u}^{2}\mathrm{d}u.
\end{equation*}%
\cite{BarndorffNielsenShephard(03test)} used this type of result to test for
jumps as this particular bipower variation is robust to jumps.
\end{example}

\section{Central limit theorem\label{sect:CLT}}

\subsection{Motivation}

It is important to be able to quantify the difference between the estimator $%
Y^{n}(g,h)$ and $Y(g,h)$. In this subsection we do this by giving a central
limit theorem for $\sqrt{n}(Y^{n}(g,h)-Y(g,h))$. We have to make some
stronger assumptions both on the process $Y$ and on the pair $(g,h)$ in
order to derive this result.

\subsection{Assumptions on the process}

We start with a variety of assumptions which strengthen (H) and (K) given in
the previous subsection.

\noindent \textbf{Assumption (H0):} We have (H) with
\begin{equation}
\sigma _{t}=\sigma _{0}+\int_{0}^{t}a_{u}^{\ast }\mathrm{d}%
u+\int_{0}^{t}\sigma _{u-}^{\ast }\mathrm{d}W_{u}+\int_{0}^{t}v_{u-}^{\ast }%
\mathrm{d}Z_{u},  \label{H'}
\end{equation}%
where $Z$ is a $d^{\prime \prime }$-dimensional L\'{e}vy process,
independent of $W$. Further, the processes $a^{\ast }$, $\sigma ^{\ast }$, $%
v^{\ast }$ are adapted c\`{a}dl\`{a}g arrays, with $a^{\ast }$ also being
predictable and locally bounded.

\noindent \textbf{Assumption (H1):} We have (H) with
\begin{eqnarray}
\sigma _{t} &=&\sigma _{0}+\int_{0}^{t}a_{u}^{\ast }\mathrm{d}%
u+\int_{0}^{t}\sigma _{u-}^{\ast }\mathrm{d}W_{u}+\int_{0}^{t}v_{u-}^{\ast }%
\mathrm{d}V_{u}  \label{assumption (V)} \\
&&+\int_{0}^{t}\int_{E}\varphi \circ w(u-,x)\left( \mu -\nu \right) \left(
\mathrm{d}u,\mathrm{d}x\right) +\int_{0}^{t}\int_{E}\left( w-\varphi \circ
w\right) \left( u-,x\right) \mu \left( \mathrm{d}u,\mathrm{d}x\right) .
\notag
\end{eqnarray}%
Here $a^{\ast }$, $\sigma ^{\ast }$, $v^{\ast }$ are adapted c\`{a}dl\`{a}g
arrays, with $a^{\ast }$ also being predictable and locally bounded. $V$ is
a $d^{\prime \prime }$-dimensional Brownian motion independent of $W$. $\mu $
is a Poisson measure on $\left( 0,\infty \right) \times E$ independent of $W$
and $V$, with intensity measure $\nu (\mathrm{d}t,\mathrm{d}x)=\mathrm{d}%
t\otimes F(\mathrm{d}x)$ and $F$ is a $\sigma $-finite measure on the Polish
space $\left( E,\mathcal{E}\right) $. $\varphi $ is a continuous truncation
function on $R^{dd^{\prime }}$ (a function with compact support, which
coincide with the identity map on the neighbourhood of $0$). Finally $%
w(\omega ,u,x)$ is a map $\Omega \times \lbrack 0,\infty )\times E$ into the
space of $d\times d^{\prime }$arrays which is $\mathcal{F}_{u}\otimes $ $%
\mathcal{E}-$measurable in $(\omega ,x)$ for all $u$ and c\`{a}dl\`{a}g in $%
u $, and such that for some sequences $\left( S_{k}\right) $ of stopping
times increasing to $+\infty $ we have%
\begin{equation*}
\sup_{\omega \in \Omega ,u<S_{k}(\omega )}\left\Vert w(\omega
,u,x)\right\Vert \leq \psi _{k}(x)\quad \text{where\quad }\int_{E}\left(
1\wedge \psi _{k}(x)^{2}\right) F(\mathrm{d}x)<\infty .
\end{equation*}

\noindent \textbf{Assumption (H2): }$\Sigma =\sigma \sigma ^{\prime }$ is
everywhere invertible.

\begin{remark}
Assumption (H1) looks quite complicated but has been setup so that the same
conditions on the coefficients can be applied both to $\sigma $ and $\Sigma
=\sigma \sigma ^{\prime }$. If there were no jumps then it would be
sufficient to employ the first line of (\ref{assumption (V)}). The
assumption (H1) is rather general from an econometric viewpoint as it allows
for flexible leverage effects, multifactor volatility effects, jumps,
non-stationarities, intraday effects, etc.
\end{remark}

\subsection{Assumptions on $g$ and $h$}

In order to derive a central limit theorem we need to impose some regularity
on $g$ and $h$.

\noindent \textbf{Assumption (K1): }$f$ is even (that is $f(x)=f(-x)$ for $%
x\in R^{d}$) and continuously differentiable, with derivatives having at
most polynomial growth.

In order to handle some of the most interesting cases of bipower variation,
where we are mostly interested in taking low powers of absolute values of
returns which may not be differentiable at zero, we sometimes need to relax
(K1). The resulting condition is quite technical and is called (K2). It is
discussed in the Appendix.

\noindent \textbf{Assumption (K2):} $f$ is even and continuously
differentiable on the complement $B^{c}$\ of a closed subset $B\subset
\mathbb{R}^{d}$ and satisfies%
\begin{equation*}
||y||\leq 1\Longrightarrow |f(x+y)-f(x)|\leq C(1+||x||^{p})||y||^{r}
\end{equation*}%
for some constants $C$, $p\geq 0$ and $r\in \left( 0,1\right] $. Moreover

a) If $r=1$ then $B$ has Lebesgue measure $0$.

b) If $r<1$ then $B$ satisfies

\begin{equation}
\left.
\begin{array}{l}
\text{for any positive definite }d\times d\text{ matrix }C\text{ and } \\
\text{any }N(0,C)\text{-random vector }U\text{ the distance }d(U,B) \\
\text{from }U\text{ to }B\text{ has a density }\psi _{C}\text{ on }R_{+},%
\text{ such that } \\
sup_{x\in R_{+},|C|+|C^{-1}|\leq A}\psi _{C}(x)<\infty \text{ for all }%
A<\infty ,%
\end{array}%
\right\}  \label{K13}
\end{equation}

\qquad\ and we have
\begin{equation}
x\in B^{c},~\Vert y\Vert \leq 1\bigwedge {\frac{d(x,B)}{2}}~~\Rightarrow
~~\left\{
\begin{array}{l}
\Vert \nabla f(x)\Vert \leq {\frac{C(1+\Vert x\Vert ^{p})}{d(x,B)^{1-r}}}, \\%
[2.5mm]
\Vert \nabla f(x+y)-\nabla f(x)\Vert \leq {\frac{C(1+\Vert x\Vert ^{p})\Vert
y\Vert }{d(x,B)^{2-r}}}.%
\end{array}%
\right.  \label{K11}
\end{equation}

\begin{remark}
These conditions accommodate the case where $f$ equals $\left\vert
x^{j}\right\vert ^{r}$: this function satisfies (K1) when $r>1$, and (K2)
when $r\in (0,1]$ (with the same $r$ of course). When $B$ is a finite union
of hyperplanes it satisfies (\ref{K13}). Also, observe that (K1) implies
(K2) with $r=1$ and $B=\emptyset $.
\end{remark}

\subsection{Central limit theorem}

Each of the following assumptions (J1) and (J2) are sufficient for the
statement of Theorem 1.3 of \cite%
{BarndorffNielsenGraversenJacodPodolskyShephard(04shiryaev)} to hold.

\noindent \textbf{Assumption (J1):} We have (H1) and $g$ and $h$ satisfy
(K1).

\noindent \textbf{Assumption (J2):} We have (H1), (H2) and $g$ and $h$
satisfy (K2).\newline

The result of the Theorem is restated in the following.

\begin{theorem}
\label{TT3}Assume at least one of (J1) and (J2) holds, then the process
\begin{equation*}
\sqrt{n}~(Y^{n}(g,h)_{t}-Y(g,h)_{t})
\end{equation*}%
converges stably in law towards a limiting process $U(g,h)$ having the form%
\begin{equation}
U(g,h)_{t}^{jk}=\sum_{j^{\prime }=1}^{d_{1}}\sum_{k^{\prime
}=1}^{d_{3}}\int_{0}^{t}\alpha (\sigma _{u},g,h)^{jk,j^{\prime }k^{\prime }}~%
\mathrm{d}B_{u}^{j^{\prime },k^{\prime }},
\end{equation}%
where%
\begin{equation*}
\sum_{l=1}^{d_{1}}\sum_{m=1}^{d_{3}}\alpha (\sigma ,g,h)^{jk,lm}\alpha
(\sigma ,g,h)^{j^{\prime }k^{\prime },lm}=A(\sigma ,g,h)^{jk,j^{\prime
}k^{\prime }},
\end{equation*}%
and%
\begin{eqnarray*}
A(\sigma ,g,h)^{jk,j^{\prime }k^{\prime }} &=&\displaystyle%
\sum_{l=1}^{d_{2}}\sum_{l^{\prime }=1}^{d_{2}}\left\{ \rho _{\sigma }\left(
g^{jl}g^{j^{\prime }l^{\prime }}\right) \rho _{\sigma }\left(
h^{lk}h^{l^{\prime }k^{\prime }}\right) +\rho _{\sigma }\left( g^{jl}\right)
\rho _{\sigma }\left( h^{l^{\prime }k^{\prime }}\right) \rho _{\sigma
}\left( g^{j^{\prime }l^{\prime }}h^{lk}\right) \right. \\
&&\displaystyle+\rho _{\sigma }\left( g^{j^{\prime }l^{\prime }}\right) \rho
_{\sigma }\left( h^{lk}\right) \rho _{\sigma }\left( g^{jl}h^{l^{\prime
}k^{\prime }}\right) \\
&&\displaystyle\left. -3\rho _{\sigma }\left( g^{jl}\right) \rho _{\sigma
}\left( g^{j^{\prime }l^{\prime }}\right) \rho _{\sigma }\left(
h^{lk}\right) \rho _{\sigma }\left( h^{l^{\prime }k^{\prime }}\right)
\right\} .
\end{eqnarray*}%
Furthermore, $B$ is a standard Wiener process which is defined on an
extension of $\left( \Omega ,\mathcal{F},\left( \mathcal{F}_{t}\right)
_{t\geq 0},P\right) $ and is independent of the $\sigma $--field $\mathcal{F}
$.
\end{theorem}

\begin{remark}
Convergence stably in law is slightly stronger than convergence in law. It
is discussed in, for example, \cite[pp. 512-518]{JacodShiryaev(03)}.
\end{remark}

\begin{remark}
Suppose $d_{3}=1$, which is the situation looked at in Example \ref{Example:
1}(e). Then $Y^{n}(g,h)_{t}$ is a vector and so the limiting law of $\sqrt{n}%
(Y^{n}(g,h)-Y(g,h))$ simplifies. It takes on the form of
\begin{equation}
U(g,h)_{t}^{j}=\sum_{j^{\prime }=1}^{d_{1}}\int_{0}^{t}\alpha (\sigma
_{u},g,h)^{j,j^{\prime }}~\mathrm{d}B_{u}^{j^{\prime }},
\end{equation}%
where%
\begin{equation*}
\sum_{l=1}^{d_{1}}\alpha (\sigma ,g,h)^{j,l}\alpha (\sigma ,g,h)^{j^{\prime
},l}=A(\sigma ,g,h)^{j,j^{\prime }}.
\end{equation*}%
Here%
\begin{eqnarray*}
A(\sigma ,g,h)^{j,j^{\prime }} &=&\displaystyle\sum_{l=1}^{d_{2}}\sum_{l^{%
\prime }=1}^{d_{2}}\left\{ \rho _{\sigma }(g^{jl}g^{j^{\prime }l^{\prime
}})\rho _{\sigma }(h^{l}h^{l^{\prime }})+\rho _{\sigma }(g^{jl})\rho
_{\sigma }(h^{l^{\prime }})\rho _{\sigma }(g^{j^{\prime }l^{\prime
}}h^{l})\right. \\
&&\displaystyle+\left. \rho _{\sigma }(g^{j^{\prime }l^{\prime }})\rho
_{\sigma }(h^{l})\rho _{\sigma }(g^{jl}h^{l^{\prime }})-3\rho _{\sigma
}(g^{jl})\rho _{\sigma }(g^{j^{\prime }l^{\prime }})\rho _{\sigma
}(h^{l})\rho _{\sigma }(h^{l^{\prime }})\right\} .
\end{eqnarray*}%
In particular, for a single point in time $t$,
\begin{equation*}
\sqrt{n}~(Y^{n}(g,h)_{t}-Y(g,h)_{t})\rightarrow MN\left(
0,\int_{0}^{t}A(\sigma _{u},g,h)\mathrm{d}u\right) ,
\end{equation*}%
where $MN$ denotes a mixed Gaussian distribution. and $A(\sigma ,g,h)$
denotes a matrix whose $j,j^{\prime }$-th element is $A(\sigma
,g,h)^{j,j^{\prime }}$.
\end{remark}

\begin{remark}
Suppose $g(y)=I$, then $A$ becomes
\begin{equation*}
A(\sigma ,g,h)^{jk,j^{\prime }k^{\prime }}=\rho _{\sigma
}(h^{jk}h^{j^{\prime }k^{\prime }})-\rho _{\sigma }(h^{jk})\rho _{\sigma
}(h^{j^{\prime }k^{\prime }}).
\end{equation*}
\end{remark}

\subsection{Leading examples of this result}

\begin{example}
Suppose $d_{1}=d_{2}=d_{3}=1$, then
\begin{equation}
U(g,h)_{t}=\int_{0}^{t}\sqrt{A(\Sigma _{u},g,h)}~\mathrm{d}B_{u},
\end{equation}%
where%
\begin{equation*}
A(\sigma ,g,h)=\rho _{\sigma }(gg)\rho _{\sigma }(hh)+2\rho _{\sigma
}(g)\rho _{\sigma }(h)\rho _{\sigma }(gh)-3\left\{ \rho _{\sigma }(g)\rho
_{\sigma }(h)\right\} ^{2}.
\end{equation*}%
We consider two concrete examples of this setup.

\noindent \textbf{(i)} Power variation. Suppose $g(y)=1$ and $%
h(y)=\left\vert y^{j}\right\vert ^{r}$ where $r>0$, then $\rho _{\sigma
}(g)=1$,
\begin{equation*}
\rho _{\sigma }(h)=\rho _{\sigma }(gh)=\mu _{r}\sigma _{j}^{r},\quad \rho
_{\sigma }(hh)=\mu _{2r}\sigma _{j}^{2r}.
\end{equation*}%
This implies that%
\begin{eqnarray*}
A(\sigma ,g,h) &=&\mu _{2r}\sigma _{j}^{2r}+2\mu _{r}^{2}\sigma
_{j}^{2r}-3\mu _{r}^{2}\sigma _{j}^{2r} \\
&=&\left( \mu _{2r}-\mu _{r}^{2}\right) \sigma _{j}^{2r} \\
&=&v_{r}\sigma _{j}^{2r},
\end{eqnarray*}%
where $v_{r}=\mathrm{Var}(\left\vert u\right\vert ^{r})$ and $u\sim N(0,1)$.
When $r=2$, this yields a central limit theorem for the realised quadratic
variation process, with
\begin{equation*}
U(g,h)_{t}=\int_{0}^{t}\sqrt{2\sigma _{j,u}^{4}}~\mathrm{d}B_{u},
\end{equation*}%
a result which appears in \cite{Jacod(94)}, \cite{MyklandZhang(05)} and,
implicitly, \cite{JacodProtter(98)}, while the case of a single value of $t$
appears in \cite{BarndorffNielsenShephard(02realised)}. For the more general
case of $r>0$ \cite{BarndorffNielsenShephard(03bernoulli)} derived, under
much stronger conditions, a central limit theorem for $U(g,h)_{1}$. Their
result ruled out leverage effects, which are allowed under Theorem \ref{TT3}%
. The finite sample behaviour of this type of limit theory is studied in,
for example, \cite{BarndorffNielsenShephard(05tom)}, \cite%
{GoncalvesMeddahi(04)} and \cite{NielsenFrederiksen(05)}.

\noindent \textbf{(ii)} Bipower variation. Suppose $g(y)=\left\vert
y^{j}\right\vert ^{r}$ and $h(y)=\left\vert y^{j}\right\vert ^{s}$ where $%
r,s>0$, then%
\begin{eqnarray*}
\rho _{\sigma }(g) &=&\mu _{r}\sigma _{j}^{r},\quad \rho _{\sigma }(h)=\mu
_{s}\sigma _{j}^{s},\quad \rho _{\sigma }(gg)=\mu _{2r}\sigma _{j}^{2r},\quad
\\
\rho _{\sigma }(hh) &=&\mu _{2s}\sigma _{j}^{2s},\quad \rho _{\sigma
}(gh)=\mu _{r+s}\sigma _{j}^{r+s}.
\end{eqnarray*}%
This implies that
\begin{eqnarray*}
A(\sigma ,g,h) &=&\mu _{2r}\sigma _{j}^{2r}\mu _{2s}\sigma _{j}^{2s}+2\mu
_{r}\sigma _{j}^{r}\mu _{s}\sigma _{j}^{s}\mu _{r+s}\sigma _{j}^{r+s}-3\mu
_{r}^{2}\sigma _{j}^{2r}\mu _{s}^{2}\sigma _{j}^{2s} \\
&=&\left( \mu _{2r}\mu _{2s}+2\mu _{r+s}\mu _{r}\mu _{s}-3\mu _{r}^{2}\mu
_{s}^{2}\right) \sigma _{j}^{2r+2s}.
\end{eqnarray*}%
In the $r=s=1$ case \cite{BarndorffNielsenShephard(03test)} derived, under
much stronger conditions, a central limit theorem for $U(g,h)_{1}$. Their
result ruled out leverage effects, which are allowed under Theorem \ref{TT3}%
. In that special case, writing
\begin{equation*}
\vartheta =\frac{\pi ^{2}}{4}+\pi -5,
\end{equation*}%
we have
\begin{equation*}
U(g,h)_{t}=\mu _{1}^{2}\int_{0}^{t}\sqrt{\left( 2+\vartheta \right) \sigma
_{j,u}^{4}}~\mathrm{d}B_{u}.
\end{equation*}
\end{example}

\begin{example}
Suppose $g=I$, $h(y)=yy^{\prime }$. Then we have to calculate%
\begin{equation*}
A(\sigma ,g,h)^{jk,j^{\prime }k^{\prime }}=\rho _{\sigma
}(h^{jk}h^{j^{\prime }k^{\prime }})-\rho _{\sigma }(h^{jk})\rho _{\sigma
}(h^{j^{\prime }k^{\prime }}).
\end{equation*}%
However,
\begin{equation*}
\rho _{\sigma }(h^{jk})=\Sigma ^{jk},\quad \rho _{\sigma
}(h^{jk}h^{j^{\prime }k^{\prime }})=\Sigma ^{jk}\Sigma ^{j^{\prime
}k^{\prime }}+\Sigma ^{jj^{\prime }}\Sigma ^{kk^{\prime }}+\Sigma
^{jk^{\prime }}\Sigma ^{kj^{\prime }},
\end{equation*}%
so
\begin{eqnarray*}
A(\sigma ,g,h)^{jk,j^{\prime }k^{\prime }} &=&\Sigma ^{jk}\Sigma ^{j^{\prime
}k^{\prime }}+\Sigma ^{jj^{\prime }}\Sigma ^{kk^{\prime }}+\Sigma
^{jk^{\prime }}\Sigma ^{kj^{\prime }}-\Sigma ^{jk}\Sigma ^{j^{\prime
}k^{\prime }} \\
&=&\Sigma ^{jj^{\prime }}\Sigma ^{kk^{\prime }}+\Sigma ^{jk^{\prime }}\Sigma
^{kj^{\prime }}.
\end{eqnarray*}%
This is the result found in \cite{BarndorffNielsenShephard(04multi)}, but
proved under stronger conditions, and is implicit in the work of \cite%
{JacodProtter(98)}.
\end{example}

\begin{example}
\label{Example:vector case}Suppose $d_{1}=d_{2}=2$, $d_{3}=1$ and $g$ is
diagonal. Then
\begin{equation}
U(g,h)_{t}^{j}=\sum_{j^{\prime }=1}^{2}\int_{0}^{t}\alpha (\sigma
_{u},g,h)^{j,j^{\prime }}~\mathrm{d}B_{u}^{j^{\prime }},
\end{equation}%
where%
\begin{equation*}
\sum_{l=1}^{2}\alpha (\sigma ,g,h)^{j,l}\alpha (\sigma ,g,h)^{j^{\prime
},l}=A(\sigma ,g,h)^{j,j^{\prime }}.
\end{equation*}%
Here%
\begin{eqnarray*}
A(\sigma ,g,h)^{j,j^{\prime }} &=&\rho _{\sigma }(g^{jj}g^{j^{\prime
}j^{\prime }})\rho _{\sigma }(h^{j}h^{j^{\prime }})+\rho _{\sigma
}(g^{jj})\rho _{\sigma }(h^{j^{\prime }})\rho _{\sigma }(g^{j^{\prime
}j^{\prime }}h^{j}) \\
&&+\rho _{\sigma }(g^{j^{\prime }j^{\prime }})\rho _{\sigma }(h^{j})\rho
_{\sigma }(g^{jj}h^{j^{\prime }})-3\rho _{\sigma }(g^{jj})\rho _{\sigma
}(g^{j^{\prime }j^{\prime }})\rho _{\sigma }(h^{j})\rho _{\sigma
}(h^{j^{\prime }}).
\end{eqnarray*}
\end{example}

\begin{example}
\label{Example:joint BPV and RV}Joint behaviour of realised QV and realised
bipower variation. This sets%
\begin{equation*}
g(y)=\left(
\begin{array}{cc}
\left\vert y^{j}\right\vert & 0 \\
0 & 1%
\end{array}%
\right) ,\quad h(y)=\left(
\begin{array}{c}
\left\vert y^{j}\right\vert \\
\left( y^{j}\right) ^{2}%
\end{array}%
\right) .
\end{equation*}%
The implication is that
\begin{equation*}
\rho _{\sigma }(g^{11})=\rho _{\sigma }(g^{22}g^{11})=\rho _{\sigma
}(g^{11}g^{22})=\mu _{1}\sigma _{j},\ \rho _{\sigma }(g^{22})=1,\ \rho
_{\sigma }(g^{11}g^{11})=\sigma _{j}^{2},\ \rho _{\sigma }(g^{22}g^{22})=1,
\end{equation*}%
\begin{equation*}
\rho _{\sigma }(h^{1})=\mu _{1}\sigma _{j},\ \rho _{\sigma }(h^{2})=\rho
_{\sigma }(h^{1}h^{1})=\sigma _{j}^{2},\ \rho _{\sigma }(h^{1}h^{2})=\rho
_{\sigma }(h^{2}h^{1})=\mu _{3}\sigma _{j}^{3},\ \rho _{\sigma
}(h^{2}h^{2})=3\sigma _{j}^{4},
\end{equation*}%
\begin{equation*}
\rho _{\sigma }(g^{11}h^{1})=\sigma _{j}^{2},\ \rho _{\sigma
}(g^{11}h^{2})=\mu _{3}\sigma _{j}^{3},\ \rho _{\sigma }(g^{22}h^{1})=\mu
_{1}\sigma _{j},\ \rho _{\sigma }(g^{22}h^{2})=\sigma _{j}^{2}.
\end{equation*}%
Thus
\begin{eqnarray*}
A(\sigma ,g,h)^{1,1} &=&\sigma _{j}^{2}\sigma _{j}^{2}+2\mu _{1}\sigma
_{j}\mu _{1}\sigma _{j}\sigma _{j}^{2}-3\mu _{1}\sigma _{j}\mu _{1}\sigma
_{j}\mu _{1}\sigma _{j}\mu _{1}\sigma _{j} \\
&=&\sigma _{j}^{4}\left( 1+2\mu _{1}^{2}-3\mu _{1}^{4}\right) =\mu
_{1}^{4}(2+\vartheta )\sigma _{j}^{4},
\end{eqnarray*}%
while%
\begin{equation*}
A(\sigma ,g,h)^{2,2}=3\sigma _{j}^{4}+2\sigma _{j}^{4}-3\sigma
_{j}^{4}=2\sigma _{j}^{4},
\end{equation*}%
and%
\begin{eqnarray*}
A(\sigma ,g,h)^{1,2} &=&\mu _{1}\sigma _{j}\mu _{3}\sigma _{j}^{3}+\mu
_{1}\sigma _{j}\sigma _{j}^{2}\mu _{1}\sigma _{j}+\mu _{1}\sigma _{j}\mu
_{3}\sigma _{j}^{3}-3\mu _{1}\sigma _{j}\mu _{1}\sigma _{j}\sigma _{j}^{2} \\
&=&2\sigma _{j}^{4}\left( \mu _{1}\mu _{3}-\mu _{1}^{2}\right) =2\mu
_{1}^{2}\sigma _{j}^{4}.
\end{eqnarray*}%
This generalises the result given in \cite{BarndorffNielsenShephard(03test)}
to the leverage case. In particular we have that
\begin{equation*}
\left(
\begin{array}{c}
U(g,h)_{t}^{1} \\
U(g,h)_{t}^{2}%
\end{array}%
\right) =\left(
\begin{array}{l}
\displaystyle\mu _{1}^{2}\int_{0}^{t}\sqrt{2\sigma _{u}^{4}}\mathrm{d}%
B_{u}^{1}+\mu _{1}^{2}\int_{0}^{t}\sqrt{\vartheta \sigma _{u}^{4}}\mathrm{d}%
B_{u}^{2} \\
\displaystyle\int_{0}^{t}\sqrt{2\sigma _{u}^{4}}\mathrm{d}B_{u}^{1}.%
\end{array}%
\right)
\end{equation*}
\end{example}

\section{Multipower variation\label{sect:multipower variation}}

A natural extension of generalised bipower variation is to generalised
multipower variation%
\begin{equation*}
Y^{n}(g)_{t}=\frac{1}{n}\sum_{i=1}^{\left\lfloor nt\right\rfloor }\left\{
\dprod\limits_{i^{\prime }=1}^{I\wedge \left( i+1\right) }g_{i^{\prime }}(%
\sqrt{n}~\Delta _{i-i^{\prime }+1}^{n}Y)\right\} .
\end{equation*}%
This measure of variation, for the $g_{i^{\prime }}$ being absolute powers,
was introduced by \cite{BarndorffNielsenShephard(03test)}.

We will be interested in studying the properties of $Y^{n}(g)_{t}$ for given
functions $\left\{ g_{i}\right\} $ with the following properties.

\noindent \textbf{Assumption (K}$^{\ast }$\textbf{):} All the $\left\{
g_{i}\right\} $ are continuous with at most polynomial growth.

The previous results suggests that if $Y$ is a Brownian semimartingale and
Assumption (K$^{\ast }$) holds then
\begin{equation*}
Y^{n}(g)_{t}\rightarrow Y(g)_{t}:=\int_{0}^{t}\dprod\limits_{i=0}^{I}\rho
_{\sigma _{u}}(g_{i})\mathrm{d}u.
\end{equation*}

\begin{example}
\textbf{(a)} Suppose $I=4$ and $g_{i}(y)=\left\vert y^{j}\right\vert $, then
$\rho _{\sigma }(g_{i})=\mu _{1}\sigma _{j}$ so
\begin{equation*}
Y(g)_{t}=\mu _{1}^{4}\int_{0}^{t}\sigma _{j,u}^{4}\mathrm{d}u,
\end{equation*}%
a scaled version of integrated quarticity. \newline
\noindent \textbf{(b)} Suppose $I=3$ and $g_{i}(y)=\left\vert
y^{j}\right\vert ^{4/3}$, then
\begin{equation*}
\rho _{\sigma }(g_{i})=\mu _{4/3}\sigma _{j}^{4/3}
\end{equation*}%
so
\begin{equation*}
Y(g)_{t}=\mu _{4/3}^{3}\int_{0}^{t}\sigma _{j,u}^{4}\mathrm{d}u.
\end{equation*}
\end{example}

\begin{example}
Of some importance is the generic case where $g_{i}(y)=\left\vert
y^{j}\right\vert ^{2/I}$, which implies
\begin{equation*}
Y(g)_{t}=\mu _{2/I}^{I}\int_{0}^{t}\sigma _{j,u}^{2}\mathrm{d}u.
\end{equation*}%
Thus this class provides an interesting alternative to realised variance as
an estimator of integrated variance. Of course it is important to know a
central limit theory for these types of quantities. \ \cite%
{BarndorffNielsenGraversenJacodPodolskyShephard(04shiryaev)} show that when
(H1) and (H2) hold then
\begin{equation*}
\sqrt{n}\left[ Y^{n}(g)_{t}-Y(g)_{t}\right] \rightarrow \int_{0}^{t}\sqrt{%
\omega _{I}^{2}\sigma _{j,u}^{4}}~\mathrm{d}B_{u},
\end{equation*}%
where%
\begin{equation*}
\omega _{I}^{2}=\mathrm{Var}\left( \dprod\limits_{i=1}^{I}\left\vert
u_{i}\right\vert ^{2/I}\right) +2\sum_{j=1}^{I-1}\mathrm{Cov}\left(
\dprod\limits_{i=1}^{I}\left\vert u_{i}\right\vert
^{2/I},\dprod\limits_{i=1}^{I}\left\vert u_{i-j}\right\vert ^{2/I}\right) ,
\end{equation*}%
with $u_{i}\sim NID(0,1)$. Clearly $\omega _{1}^{2}=2$, while recalling that
$\mu _{1}=\sqrt{2/\pi }$,
\begin{eqnarray*}
\omega _{2}^{2} &=&\mathrm{Var}(\left\vert u_{1}\right\vert \left\vert
u_{2}\right\vert )+2\mathrm{Cov}(\left\vert u_{1}\right\vert \left\vert
u_{2}\right\vert ,\left\vert u_{2}\right\vert \left\vert u_{3}\right\vert )
\\
&=&1+2\mu _{1}^{2}-3\mu _{1}^{4},
\end{eqnarray*}

\noindent and
\begin{eqnarray*}
\omega _{3}^{2} &=&\mathrm{Var}(\left( \left\vert u_{1}\right\vert
\left\vert u_{2}\right\vert \left\vert u_{3}\right\vert \right) ^{2/3})+2%
\mathrm{Cov}(\left( \left\vert u_{1}\right\vert \left\vert u_{2}\right\vert
\left\vert u_{3}\right\vert \right) ^{2/3},\left( \left\vert
u_{2}\right\vert \left\vert u_{3}\right\vert \left\vert u_{4}\right\vert
\right) ^{2/3}) \\
&&+2\mathrm{Cov}(\left( \left\vert u_{1}\right\vert \left\vert
u_{2}\right\vert \left\vert u_{3}\right\vert \right) ^{2/3},\left(
\left\vert u_{3}\right\vert \left\vert u_{4}\right\vert \left\vert
u_{5}\right\vert \right) ^{2/3}) \\
&=&\left( \mu _{4/3}^{3}-\mu _{2/3}^{6}\right) +2\left( \mu _{4/3}^{2}\mu
_{2/3}^{2}-\mu _{2/3}^{6}\right) +2\left( \mu _{4/3}\mu _{2/3}^{4}-\mu
_{2/3}^{6}\right) .
\end{eqnarray*}
\end{example}

\begin{example}
The law of large numbers and the central limit theorem also hold for linear
combinations of processes like $Y(g)$ above. For example one may denote by $%
\zeta^n_i$ the $d\times d$ matrix whose $(k,l)$ entry is $%
\sum_{j=0}^{d-1}\Delta^n_{i+j}Y^k\Delta^n_{i+j}Y^l$. Then
\begin{equation*}
Z^n_t=\frac{n^{d-1}}{d!}\sum_{i=1}^{[nt]}\det(\zeta^n_i)
\end{equation*}
is a linear combinations of processes $Y^n(g)$ for functions $g_l$ being of
the form $g_l(y)=y^jy^k$. It is proved in \cite{JacodLejayTalay(05)} that
under (H)
\begin{equation*}
Z^n_t\rightarrow Z_t:=\int_0^t \det(\sigma_u\sigma^{\prime}_u)du
\end{equation*}
in probability, whereas under (H1) and (H2) the associated CLT is the
following convergence in law:
\begin{equation*}
\sqrt{n}(Z^n_t-Z_t)\rightarrow \int_0^t\sqrt{\Gamma(\sigma_u)}~dB_u,
\end{equation*}
where $\Gamma(\sigma)$ denotes the covariance of the variable $%
\det(\zeta)/d! $, and $\zeta$ is a $d\times d$ matrix whose $(k,l)$ entry is
$\sum_{j=0}^{d-1}U_j^kU_j^l$ and the $U_j$'s are i.i.d. centered Gaussian
vectors with covariance $\sigma\sigma^{\prime}$.

This kind of result may be used for testing whether the rank of the
diffusion coefficient is everywhere smaller than $d$ (in which case one
could use a model with a $d^{\prime}<d$ for the dimension of the driving
Wiener process $W$).
\end{example}

\section{Conclusion}

This paper provides some rather general limit results for realised
generalised bipower variation. In the case of power variation and bipower
variation the results are proved under much weaker assumptions than those
which have previously appeared in the literature. In particular the
no-leverage assumption is removed, which is important in the application of
these results to stock data.

There are a number of open questions. It is rather unclear how
econometricians might exploit the generality of the $g$ and $h$ functions to
learn about interesting features of the variation of price processes. It
would be interesting to know what properties $g$ and $h$ must possess in
order for these statistics to be robust to finite activity and infinite
activity jumps. A challenging extension is to construct a version of
realised generalised bipower variation which is robust to market
microstructure effects. Following the work on the realised volatility there
are two leading strategies which may be able to help: the kernel based
approach, studied in detailed by \cite%
{BarndorffNielsenHansenLundeShephard(04)}, and the subsampling approach of
\cite{ZhangMyklandAitSahalia(03)} and \cite{Zhang(04)}. In the realised
volatility case these methods are basically equivalent, however it is
perhaps the case that the subsampling method is easier to extend to the
non-quadratic case.

\section{Acknowledgments}

Ole E. Barndorff-Nielsen's work is supported by CAF (\texttt{www.caf.dk}),
which is funded by the Danish Social Science Research Council. Neil
Shephard's research is supported by the UK's ESRC through the grant
\textquotedblleft High frequency financial econometrics based upon power
variation.\textquotedblright\

\section{Proof of Theorem \protect\ref{TT3}}

\subsection{Strategy for the proof}

Below we give a fairly detailed account of the basic techniques in the proof
of Theorem \ref{TT3}, in the one-dimensional case and under some relatively
minor simplifying assumptions. Throughout we set $h=1$ for the main
difficulty in the proof is being able to deal with the generality in the $g$
function. Once that has been mastered the extension to the bipower measure
is not a large obstacle. We refer the reader to \cite%
{BarndorffNielsenGraversenJacodPodolskyShephard(04shiryaev)} for readers who
wish to see the more general case. In this subsection we provide a brief
outline of the content of the Section.

The aim of this Section is to show that%
\begin{equation}
\sqrt{n}\left( \frac{1}{n}\,\sum_{i=1}^{[nt]}g\left( \sqrt{n}\,\triangle
_{i}^{n}Y\right) -\int_{0}^{t}\rho _{\sigma _{u}}(g)\right) \rightarrow
\int_{0}^{t}\sqrt{\rho _{\sigma _{u}}(g^{2})-\rho _{\sigma _{u}}(g)^{2}}\;%
\mathrm{d}B_{u}  \label{eqn 0}
\end{equation}%
where $B$\ is a Brownian motion independent of the process $Y$\ and the
convergence is (stably) in law. This case is important for the extension to
realised generalised bipower (and multipower) variation is relatively simple
once this fundamental result is established.

The proof of this result is done\ in a number of steps, some of them
following fairly standard reasoning, others requiring special techniques.

The first step is to rewrite the left hand side of (\ref{eqn 0}) as follows%
\begin{eqnarray*}
&&\sqrt{n}\left( \frac{1}{n}\,\sum_{i=1}^{[nt]}g(\sqrt{n}\,\triangle
_{i}^{n}Y)-\int_{0}^{t}\rho _{\sigma _{u}}(g)\mathrm{d}u\right) \\
&=&\frac{1}{\sqrt{n}}\,\sum_{i=1}^{[nt]}\left\{ g(\sqrt{n}\,\triangle
_{i}^{n}Y)-\mathrm{E}\left[ g(\triangle _{i}^{n}Y)\,|\,\mathcal{F}_{\frac{i-1%
}{n}}\right] \right\} \, \\
&&+\sqrt{n}\left( \frac{1}{n}\,\sum_{i=1}^{[nt]}\mathrm{E}\left[ g(\triangle
_{i}^{n}Y)\,|\,\mathcal{F}_{\frac{i-1}{n}}\right] \,-\int_{0}^{t}\rho
_{\sigma _{u}}(g)\mathrm{d}u\right) .
\end{eqnarray*}%
It is rather straightforward to show that the first term of the right hand
side satisfies
\begin{equation*}
\frac{1}{\sqrt{n}}\,\sum_{i=1}^{[nt]}\left\{ \,g(\sqrt{n}\,\triangle
_{i}^{n}Y)-\mathrm{E}\left[ g(\triangle _{i}^{n}Y)\,|\,\mathcal{F}_{\frac{i-1%
}{n}}\right] \right\} \rightarrow \int_{0}^{t}\sqrt{\rho _{\sigma
_{u}}(g^{2})-\rho _{\sigma _{u}}(g)^{2}}\mathrm{d}B_{u}.
\end{equation*}%
Hence what remains is to verify that%
\begin{equation}
\sqrt{n}\left( \frac{1}{n}\,\sum_{i=1}^{[nt]}\mathrm{E}\left[ g(\triangle
_{i}^{n}Y)\,|\,\mathcal{F}_{\frac{i-1}{n}}\right] \,-\int_{0}^{t}\rho
_{\sigma _{u}}(g)\mathrm{d}u\right) \rightarrow 0.  \label{2}
\end{equation}%
We have%
\begin{eqnarray}
&&\sqrt{n}\left( \frac{1}{n}\,\sum_{i=1}^{[nt]}\mathrm{E}\left[ g(\triangle
_{i}^{n}Y)\,|\,\mathcal{F}_{\frac{i-1}{n}}\right] \,-\int_{0}^{t}\rho
_{\sigma _{u}}(g)\mathrm{d}u\right)  \notag \\
&=&\frac{1}{\sqrt{n}}\,\sum_{i=1}^{[nt]}\mathrm{E}\left[ g(\triangle
_{i}^{n}Y)\,|\,\mathcal{F}_{\frac{i-1}{n}}\right] \,-\sqrt{n}%
\sum_{i=1}^{[nt]}\int_{(i-1)/n}^{i/n}\rho _{\sigma _{u}}(g)\mathrm{d}u
\notag \\
&&+\sqrt{n}\left( \sum_{i=1}^{[nt]}\int_{(i-1)/n}^{i/n}\rho _{\sigma _{u}}(g)%
\mathrm{d}u-\int_{0}^{t}\rho _{\sigma _{u}}(g)\mathrm{d}u\right)  \label{3}
\end{eqnarray}%
where%
\begin{equation*}
\sqrt{n}\left\{ \sum_{i=1}^{[nt]}\int_{(i-1)/n}^{i/n}\rho _{\sigma _{u}}(g)%
\mathrm{d}u-\int_{0}^{t}\rho _{\sigma _{u}}(g)\mathrm{d}u\right\}
\rightarrow 0.
\end{equation*}%
The first term on the right hand side of (\ref{3}) is now split into the
difference of%
\begin{equation}
\frac{1}{\sqrt{n}}\,\sum_{i=1}^{[nt]}\left\{ \mathrm{E}\left[ g(\triangle
_{i}^{n}Y)\,|\,\mathcal{F}_{\frac{i-1}{n}}\right] \,-\rho _{\frac{i-1}{n}%
}\right\}  \label{4}
\end{equation}%
where
\begin{equation*}
\rho _{\frac{i-1}{n}}=\rho _{\sigma _{\frac{i-1}{n}}}(g)=\mathrm{E}\left[
g(\sigma _{\frac{i-1}{n}}\triangle _{i}^{n}W)\,|\,\mathcal{F}_{\frac{i-1}{n}}%
\right]
\end{equation*}%
and%
\begin{equation}
\sqrt{n}\sum_{i=1}^{[nt]}\int_{(i-1)/n}^{i/n}\left\{ \rho _{\sigma _{u}}(g)%
\mathrm{d}u-\rho _{\frac{i-1}{n}}\right\} \mathrm{d}u.  \label{4b}
\end{equation}%
It is rather easy to show that (\ref{4}) tends to $0$ in probability
uniformly in $t$. The challenge is thus to show the same result holds for (%
\ref{4b}).

To handle (\ref{4b}) one splits the individual terms in the sum into%
\begin{equation}
\sqrt{n}\ \Phi ^{\prime }\left( \sigma _{\frac{i-1}{n}}\right)
\int_{(i-1)/n}^{i/n}\left( \sigma _{u}-\sigma _{\frac{i-1}{n}}\right)
\mathrm{\,d}u  \label{5}
\end{equation}%
plus%
\begin{equation}
\sqrt{n}\,\int_{(i-1)/n}^{i/n}\,\left\{ \Phi (\sigma _{u})-\Phi \left(
\sigma _{\frac{i-1}{n}}\right) -\Phi ^{\prime }\left( \sigma _{\frac{i-1}{n}%
}\right) \cdot \left( \sigma _{u}-\sigma _{\frac{i-1}{n}}\right) \right\}
\,\,\mathrm{d}u,  \label{6}
\end{equation}%
where $\Phi (x)$ is a shorthand for $\rho _{x}(g)$ and $\Phi ^{\prime }(x)$
denotes the derivative with respect to $x$.\ That (\ref{6})\ tends to $0$\
may be shown via splitting it into two terms, each of which tends to $0$\ as
is verified by a sequence of inequalities, using in particular Doob's
inequality. To prove that (\ref{5}) converges to $0$, again one splits, this
time into three terms, using the differentiability of $g$\ in the relevant
regions and the mean value theorem for differentiable functions. The two
first of these terms can be handled by relatively simple means, the third
poses the most difficult part of the whole proof and is treated via
splitting it into seven parts. It is at this stage that the assumption that $%
g$\ be even comes into play and is crucial.

This section has six other subsections. In subsection \ref%
{subsection:conventions} we introduce our basic notation, while in \ref%
{subsection:model assumptions} we set out the model and review the
assumptions we use. In subsection \ref{subsection:main result} we state the
theorem we will prove. Subsections \ref{subsection: intermediate limiting
results}, \ref{subsection:13b} and \ref{subsection:proof of 13a} give the
proofs of the successive steps.

\subsection{Notational conventions \label{subsection:conventions}}

All processes mentioned in the following are defined on a given filtered
probability space $(\Omega ,\mathcal{F},(\mathcal{F}_{t}),P)$. We shall in
general use standard notation and conventions. For instance, given a process
$(Z_{t})$ we write
\begin{equation*}
\triangle _{i}^{n}Z:=Z_{\frac{i}{n}}-Z_{\frac{i-1}{n}},\ \ \ i,n\geq 1.
\end{equation*}

We are mainly interested in convergence in law of sequences of c\`{a}dl\`{a}%
g processes. In fact all results to be proved will imply convergence `stably
in law' which is a slightly stronger notion. For this we shall use the
notation
\begin{equation*}
(Z_{t}^{n})\rightarrow (Z_{t}),
\end{equation*}%
where $(Z_{t}^{n})$ and $(Z_{t})$ are given c\`{a}dl\`{a}g processes.
Furthermore we shall write
\begin{equation*}
(Z_{t}^{n})\overset{P}{\rightarrow }0\ \ \ \text{meaning}\ \ \sup_{0\leq
s\leq t}|Z_{s}^{n}|\rightarrow 0\ \ \mbox{\rm in\ probability\ for\ all}\
t\geq 0,
\end{equation*}%
\begin{equation*}
(Z_{t}^{n})\overset{P}{\rightarrow }(Z_{t})\ \ \ \text{meaning}\ \
(Z_{t}^{n}-Z_{t})\overset{P}{\rightarrow }0.
\end{equation*}%
Often
\begin{equation*}
Z_{t}^{n}=\sum_{i=1}^{[nt]}a_{i}^{n}\ \ \ \text{for all}\ t\geq 0,
\end{equation*}%
where the $a_{i}^{n}$'s are $\mathcal{F}_{\frac{i-1}{n}}$-measurable. Recall
here that given c\`{a}dl\`{a}g processes $(Z_{t}^{n}),\,(Y_{t}^{n})$ and $%
(Z_{t})$ we have
\begin{equation*}
(Z_{t}^{n})\rightarrow (Z_{t})\ \ \text{if}\ \ (Z_{t}^{n}-Y_{t}^{n})\overset{%
P}{\rightarrow }0\ \ \text{and}\ \ (Y_{t}^{n})\rightarrow (Z_{t}).\vspace{1mm%
}
\end{equation*}

Moreover, for $h:\mathbf{R}\rightarrow \mathbf{R}$ Borel measurable of at
most polynomial growth we note that $x\mapsto \rho _{x}(h)$ is locally
bounded and continuous if $h$ is continuous at $0$.\vspace{1mm}\newline
In what follows many arguments will consist of a series of estimates of
terms indexed by $i,n$ and $t$. In these estimates we shall denote by $C$ a
finite constant which may vary from place to place. Its value will depend on
the constants and quantities appearing in the assumptions of the model but
it is always independent of $i,n$ and $t$.

\subsection{Model and basic assumptions \label{subsection:model assumptions}}

Throughout the following $(W_{t})$ denotes a $((\mathcal{F}_{t}),P)$-Wiener
process and $(\sigma _{t})$ a given c\`{a}dl\`{a}g $(\mathcal{F}_{t})$%
-adapted process. Define
\begin{equation*}
Y_{t}:=\int_{0}^{t}\sigma _{s-}\,\mathrm{d}W_{s}\ \ \ \ \ t\geq 0,
\end{equation*}%
implying that is $(Y_{t})$ is a continuous local martingale. We have deleted
the drift of the $\left( Y_{t}\right) $ process as taking care of it is a
simple technical task, while its presence increase the clutter of the
notation. Our aim is to study the asymptotic behaviour of the processes
\begin{equation*}
\{(X_{t}^{n}(g))\,|\,n\geq 1\,\}
\end{equation*}%
where%
\begin{equation*}
X_{t}^{n}(g)=\frac{1}{n}\,\sum_{i=1}^{[nt]}g(\sqrt{n}\,\triangle
_{i}^{n}Y),\ \ \ t\geq 0,\,n\geq 1.
\end{equation*}%
Here $g:\mathbf{R}\rightarrow \mathbf{R}$ is a given continuous function of
at most polynomial growth. We are especially interested in $g$'s of the form
$x\mapsto |x|^{r}\ (r>0)$ but we shall keep the general notation since
nothing is gained in simplicity by assuming that $g$ is of power form. We
shall throughout the following assume that $g$ furthermore satisfies the
following.

\begin{assumption}
\textbf{(K)}: $g$ is an even function and continuously differentiable in $%
B^{c}$ where $B\subseteq \mathbf{R}$ is a closed Lebesgue null-set and $%
\exists \ M,\,p\geq 1$ such that
\begin{equation*}
|g(x+y)-g(x)|\leq M(1+|x|^{p}+|y|^{p})\cdot |y|\ ,
\end{equation*}%
for all $x,y\in \mathbf{R}$.
\end{assumption}

\begin{remark}
The assumption (K) implies, in particular, that if $x\in B^{c}$ then
\begin{equation*}
|g^{\prime }(x)|\leq M(1+|x|^{p}).\vspace{1mm}
\end{equation*}%
Observe that only power functions corresponding to $r\geq 1$ do satisfy (K).
The remaining case $0<r<1$ requires special arguments which will be omitted
here\vspace{1mm} (for details see \cite%
{BarndorffNielsenGraversenJacodPodolskyShephard(04shiryaev)}).
\end{remark}

In order to prove the CLT-theorem we need some additional structure on the
volatility process $(\sigma _{t})$. A natural set of assumptions would be
the following.

\begin{assumption}
\textbf{(H0)}: $(\sigma _{t})$ can be written as
\begin{equation*}
\sigma _{t}=\sigma _{0}+\int_{0}^{t}a_{s}^{\ast }\,\mathrm{d}%
s+\int_{0}^{t}\sigma _{s}^{\ast }\,\mathrm{d}W_{s}+\int_{0}^{t}v_{s-}^{\ast
}\,\mathrm{d}Z_{s}
\end{equation*}%
where $(Z_{t})$ is a $((\mathcal{F}_{t}),P)$-L\'{e}vy process independent of
$(W_{t})$ and $(\sigma _{t}^{\ast })$ and $(v_{t}^{\ast })$ are adapted c%
\`{a}dl\`{a}g processes and $(a_{t}^{\ast })$ a predictable locally bounded
process.
\end{assumption}

However, in modelling volatility it is often more natural to define $(\sigma
_{t}^{2})$ as being of the above form, i.e.%
\begin{equation*}
\sigma _{t}^{2}=\sigma _{0}^{2}+\int_{0}^{t}a_{s}^{\ast }\,\mathrm{d}%
s+\int_{0}^{t}\sigma _{s}^{\ast }\,\mathrm{d}W_{s}+\int_{0}^{t}v_{s-}^{\ast
}\,\mathrm{d}Z_{s}.
\end{equation*}%
Now this does not in general imply that $(\sigma _{t})$ has the same form;
therefore we shall replace (H0) by the more general structure given by the
following assumption.

\begin{assumption}
\textbf{(H1)}: $(\sigma _{t})$ can be written, for $t\geq 0$, as%
\begin{equation*}
\begin{array}{lll}
\sigma _{t} & = & \displaystyle\sigma _{0}+\int_{0}^{t}a_{s}^{\ast }\,%
\mathrm{d}s+\int_{0}^{t}\sigma _{s}^{\ast }\,\mathrm{d}W_{s}+%
\int_{0}^{t}v_{s-}^{\ast }\,\mathrm{d}V_{s} \\
&  & \displaystyle+\int_{0}^{t}\int_{E}q\circ \phi (s-,x)\,(\mu -\nu )(%
\mathrm{d}s\,\mathrm{d}x) \\
&  & \displaystyle+\int_{0}^{t}\int_{E}\ \left\{ \phi (s-,x)-q\circ \phi
(s-,x)\right\} \,\mu (\mathrm{d}s\,\mathrm{d}x).%
\end{array}%
\end{equation*}%
Here $(a_{t}^{\ast }),\,(\sigma _{t}^{\ast })$ and $(v_{t}^{\ast })$ are as
in (H0) and $(V_{t})$ is another $((\mathcal{F}_{t}),P)$-Wiener process
independent of $(W_{t})$ while $q$ is a continuous truncation function on $%
\mathbf{R}$, i.e. a function with compact support coinciding with the
identity on a neighbourhood of $0$. Further $\mu $ is a Poisson random
measure on $(0,\infty )\times E$ independent of $(W_{t})$ and $(V_{t})$ with
intensity measure $\nu (\mathrm{d}s\,\mathrm{d}x)=\mathrm{d}s\otimes F(%
\mathrm{d}x)$, $F$ being a $\sigma $-finite mea\-sure on a measurable space $%
(E,\mathcal{E})$ and
\begin{equation*}
(\omega ,s,x)\mapsto \phi (\omega ,s,x)
\end{equation*}%
is a map from $\Omega \times \,[\,0,\infty )\times E$ into $\mathbf{R}$
which is $\mathcal{F}_{s}\otimes \mathcal{E}$ measurable in $(\omega ,x)$
for all $s$ and c\`{a}dl\`{a}g in $s$, satisfying furthermore that for some
sequence of stopping times $(S_{k})$ increasing to $+\infty $ we have for
all $k\geq 1$
\begin{equation*}
\int_{E}\left\{ 1\wedge \psi _{k}(x)^{2}\right\} \,F(\mathrm{d}x)<\infty ,
\end{equation*}%
where%
\begin{equation*}
\psi _{k}(x)=\sup_{\omega \in \Omega ,\,s<S_{k}(\omega )}|\phi (\omega
,s,x)|.
\end{equation*}
\end{assumption}

\begin{remark}
(H1) is weaker than (H0), and if $(\sigma _{t}^{2})$ satisfies (H1) then so
does $(\sigma _{t})$.\newline
\end{remark}

Finally we shall also assume a non-degeneracy in the model.

\begin{assumption}
\textbf{(H2)}: $(\sigma _{t})$ satisfies
\begin{equation*}
0<\sigma _{t}^{2}(\omega )\ \text{for all}\ (t,\omega ).\vspace{1mm}
\end{equation*}
\end{assumption}

According to general stochastic analysis theory it is known that to prove
convergence in law of a sequence $(Z_{t}^{n})$ of c\`{a}dl\`{a}g processes
it suffices to prove the convergence of each of the stopped processes $%
(Z_{T_{k}\wedge t}^{n})$ for at least one sequence of stopping times $%
(T_{k}) $ increasing to $+\infty $. Applying this together with standard
localisation techniques (for details see \cite%
{BarndorffNielsenGraversenJacodPodolskyShephard(04shiryaev)}), we may assume
that the following more restrictive assumptions are satisfied.

\begin{assumption}
\textbf{(H1a)}: $(\sigma _{t})$ can be written as
\begin{equation*}
\sigma _{t}=\sigma _{0}+\int_{0}^{t}a_{s}^{\ast }\,ds+\int_{0}^{t}\sigma
_{s-}^{\ast }\,\mathrm{d}W_{s}+\int_{0}^{t}v_{s-}^{\ast }\,\mathrm{d}%
V_{s}+\int_{0}^{t}\int_{E}\phi (s-,x)(\mu -\nu )(\mathrm{d}s\,\mathrm{d}x)\
\ \ t\geq 0.\vspace{1mm}
\end{equation*}%
Here $(a_{t}^{\ast }),\,(\sigma _{t}^{\ast })$ and $(v_{t}^{\ast })$ are
real valued uniformly bounded c\`{a}dl\`{a}g $(\mathcal{F}_{t})$-adapted
proces\-ses; $(V_{t})$ is another $((\mathcal{F}_{t}),P)$-Wiener process
independent of $(W_{t})$. Further $\mu $ is a Poisson random measure on $%
(0,\infty )\times E$ independent of $(W_{t})$ and $(V_{t})$ with intensity
measure $\nu (\mathrm{d}s\,\mathrm{d}x)=\mathrm{d}s\otimes F(\mathrm{d}x)$, $%
F$ being a $\sigma $-finite mea\-sure on a measurable space $(E,\mathcal{E})$
and
\begin{equation*}
(\omega ,s,x)\mapsto \phi (\omega ,s,x)
\end{equation*}%
is a map from $\Omega \times \,[\,0,\infty )\times E$ into $\mathbf{R}$
which is $\mathcal{F}_{s}\otimes \mathcal{E}$ measurable in $(\omega ,x)$
for all $s$ and c\`{a}dl\`{a}g in $s$, satisfying furthermore
\begin{equation*}
\psi (x)=\sup_{\omega \in \Omega ,\,s\geq 0}|\phi (\omega ,s,x)|\leq
M<\infty \ \ \text{and}\ \ \int \psi (x)^{2}\,F(\mathrm{d}x)<\infty .\vspace{%
1mm}
\end{equation*}
\end{assumption}

Likewise, by a localisation argument, we may assume

\begin{assumption}
\textbf{(H2a)}: $(\sigma _{t})$ satisfies
\begin{equation*}
a<\sigma _{t}^{2}(\omega )<b\ \ \ \text{for all}\ (t,\omega )\ \text{for some%
}\ a,b\in (0,\infty ).\vspace{1mm}
\end{equation*}
\end{assumption}

Observe that under the more restricted assumptions $(Y_{t})$ is a continuous
martingale having moments of all orders and $(\sigma _{t})$ is represented
as a sum of three square integrable martingales plus a continuous process of
bounded variation. Furthermore, the increments of the increasing processes
corresponding to the three martingales and of the bounded variation process
are dominated by a constant times $\triangle t$, implying in particular that
\begin{equation}
\mathrm{E}\left[ \,\left\vert \sigma _{v}-\sigma _{u}\right\vert ^{2}\right]
\leq C\,(v-u),\ \ \ \ \text{for all}\ 0\leq u<v.\vspace{2mm}  \label{8}
\end{equation}

\subsection{Main result \label{subsection:main result}\newline
}

As already mentioned, our aim is to show the following special version of
the general CLT-result given as Theorem \ref{TT3}.

\begin{theorem}
\vspace{2mm}Under assumptions (K), (H1a) and (H2a), there exists a Wiener
process $(B_{t})$ defined on some extension of $(\Omega ,\mathcal{F},(%
\mathcal{F}_{t}),P)$ and independent of $\mathcal{F}$ such that%
\begin{equation}
\left( \sqrt{n}\left( \,\ \frac{1}{n}\,\sum_{i=1}^{[nt]}g(\sqrt{n}%
\,\triangle _{i}^{n}Y)-\int_{0}^{t}\rho _{\sigma _{u}}(g)\,\mathrm{d}%
u\,\right) \right) \rightarrow \int_{0}^{t}\sqrt{\rho _{\sigma
_{u-}}(g^{2})-\rho _{\sigma _{u-}}(g)^{2}}\,\mathrm{d}B_{u}.
\label{Main result}
\end{equation}%
\emph{\ }
\end{theorem}

Introducing the notation%
\begin{equation*}
U_{t}(g)=\int_{0}^{t}\sqrt{\rho _{\sigma _{u-}}(g^{2})-\rho _{\sigma
_{u-}}(g)^{2}}\,\mathrm{d}B_{u}\ \ \ t\geq 0\vspace{1mm}
\end{equation*}%
we may reexpress (\ref{Main result}) as
\begin{equation}
\left( \sqrt{n}\,\left( X_{t}^{n}(g)-\int_{0}^{t}\sigma _{u}(g)\,\mathrm{d}%
u\right) \,\right) \rightarrow (U_{t}(g)).  \label{Main result reform}
\end{equation}%
To prove this, introduce the set of variables $\{\beta
_{i}^{n}\,|\,i,\,n\geq 1\}$ given by
\begin{equation*}
\beta _{i}^{n}=\sqrt{n}\cdot \sigma _{\frac{i-1}{n}}\cdot \triangle
_{i}^{n}W,\ \ \ i,\,n\geq 1.
\end{equation*}

The $\beta _{i}^{n}$'s should be seen as approximations to $\sqrt{n}%
\,\triangle _{i}^{n}Y$. In fact, since
\begin{equation*}
\sqrt{n}\,\triangle _{i}^{n}Y-\beta _{i}^{n}=\sqrt{n}\,\int_{(i-1)/n}^{i/n}(%
\sigma _{s}-\sigma _{\frac{i-1}{n}})\,\mathrm{d}W_{s}
\end{equation*}%
and $(\sigma _{t})$ is uniformly bounded, a straightforward application of (%
\ref{8}) and the Burkholder-Davis-Gundy-inequalities (e.g. \cite[pp. 160-171]%
{RevuzYor(99)}) gives for every $p>0$ the following simple estimates.
\begin{equation}
\mathrm{E}\left[ \,|\sqrt{n}\,\triangle _{i}^{n}Y-\beta _{i}^{n}|^{p}\,|\,%
\mathcal{F}_{\frac{i-1}{n}}\right] \leq \frac{C_{p}}{n^{p\wedge 1}}
\end{equation}%
and
\begin{equation}
\mathrm{E}\left[ \,|\sqrt{n}\,\triangle _{i}^{n}Y|^{p}+|\beta
_{i}^{n}|^{p}\,|\,\mathcal{F}_{\frac{i-1}{n}}\right] \leq C_{p}\vspace{1mm}
\label{12}
\end{equation}%
for all $i,n\geq 1$. Observe furthermore that
\begin{equation*}
\mathrm{E}\left[ g(\beta _{i}^{n})\,|\,\mathcal{F}_{\frac{i-1}{n}}\right]
=\rho _{\sigma _{\frac{i-1}{n}}}(g),\ \ \ \text{for all}\ i,\,n\geq 1.%
\vspace{1mm}
\end{equation*}

Introduce for convenience, for each $t>0$ and $n\geq 1$, the shorthand
notation
\begin{equation*}
U_{t}^{n}(g)=\frac{1}{\sqrt{n}}\,\sum_{i=1}^{[nt]}\,\,\left\{ g(\sqrt{n}%
\,\triangle _{i}^{n}Y)-\mathrm{E}\left[ g(\sqrt{n}\,\triangle _{i}^{n}Y)\,|\,%
\mathcal{F}_{\frac{i-1}{n}}\right] \right\} \,
\end{equation*}%
and
\begin{equation*}
\tilde{U}_{t}^{n}(g)=\frac{1}{\sqrt{n}}\,\sum_{i=1}^{[nt]}\,\,\left\{
g(\beta _{i}^{n})-\rho _{\sigma _{\frac{i-1}{n}}}(g)\right\} =\frac{1}{\sqrt{%
n}}\,\sum_{i=1}^{[nt]}\,\,\left\{ g(\beta _{i}^{n})-\mathrm{E}\left[ g(\beta
_{i}^{n})\,|\,\mathcal{F}_{\frac{i-1}{n}}\right] \right\} \,.\vspace{1mm}
\end{equation*}%
The asymptotic behaviour of $(\tilde{U}_{t}^{n}(g))$ is well known. More
precisely under the the given assumptions\thinspace (\thinspace in fact much
less is needed\thinspace ) we have
\begin{equation*}
(U_{t}^{n}(g))\rightarrow (U_{t}(g)).\vspace{1mm}
\end{equation*}%
This result is a rather straightforward consequence of \cite[Theorem IX.7.28]%
{JacodShiryaev(03)}. Thus, if $(U_{t}^{n}(g)-\tilde{U}_{t}^{n}(g))\overset{P}%
{\rightarrow }0$ we may deduce the following result.

\begin{theorem}
\label{theorem B}\ \ \emph{Let $(B_{t})$ and $(U_{t}(g))$ be as above. Then}
\begin{equation*}
(\tilde{U}_{t}^{n}(g))\rightarrow (U_{t}(g)).\vspace{1mm}
\end{equation*}
\end{theorem}

\noindent \textbf{Proof.}

As pointed out just above it is enough to prove that
\begin{equation*}
(U_{t}^{n}(g)-\tilde{U}_{t}^{n}(g))\overset{P}{\rightarrow }0.
\end{equation*}%
But for $t\geq 0$ and $n\geq 1$
\begin{equation*}
U_{t}^{n}(g)-\tilde{U}_{t}^{n}(g)=\sum_{i=1}^{[nt]}\,\left( \xi _{i}^{n}-%
\mathrm{E}\left[ \xi _{i}^{n}\,|\,\mathcal{F}_{\frac{i-1}{n}}\right] \right)
\end{equation*}%
where
\begin{equation*}
\xi _{i}^{n}=\frac{1}{\sqrt{n}}\left\{ g(\sqrt{n}\triangle
_{i}^{n}Y)-g(\beta _{i}^{n})\right\} ,\ \ \ i,n\geq 1.
\end{equation*}%
Thus we have to prove
\begin{equation*}
\left( \,\sum_{i=1}^{[nt]}\,\left\{ \xi _{i}^{n}-\mathrm{E}\left[ \xi
_{i}^{n}\,|\,\mathcal{F}_{\frac{i-1}{n}}\right] \right\} \right) \overset{P}{%
\rightarrow }0.
\end{equation*}%
But, as the left hand side of this relation is a sum of martingale
differences, this is implied by Doob's inequality (e.g. \cite[pp. 54-55]%
{RevuzYor(99)}) if for all $t>0$
\begin{equation*}
\sum_{i=1}^{[nt]}\,\mathrm{E}[(\xi _{i}^{n})^{2}]=\mathrm{E}%
[\,\sum_{i=1}^{[nt]}\,\mathrm{E}[(\xi _{i}^{n})^{2}\,|\,\mathcal{F}_{\frac{%
i-1}{n}}]\,]\rightarrow 0\ \ \ \text{as}\ n\rightarrow \infty .\vspace{1mm}
\end{equation*}%
Fix $t>0$. Using the Cauchy-Schwarz inequality and the
Burkholder-Davis-Gundy inequalities we have for all $i,n\geq 1$.
\begin{eqnarray*}
\mathrm{E}\left[ (\xi _{i}^{n})^{2}\,|\,\mathcal{F}_{\frac{i-1}{n}}\right]
&=&\frac{1}{n}\,\mathrm{E}\left[ \left\{ g(\sqrt{n}\triangle
_{i}^{n}Y)-\beta _{i}^{n}+\beta _{i}^{n}-g(\beta _{i}^{n})\right\} ^{2}\,|\,%
\mathcal{F}_{\frac{i-1}{n}}\right] \\
&\leq &\frac{C}{n}\,\mathrm{E}\left[ \,(1+|\sqrt{n}\triangle
_{i}^{n}Y|^{p}+|\beta _{i}^{n}|^{p})^{2}\cdot (\sqrt{n}\triangle
_{i}^{n}Y-\beta _{i}^{n})^{2}\,|\,\mathcal{F}_{\frac{i-1}{n}}\right] \\
&\leq &\frac{C}{n}\,\sqrt{\mathrm{E}\left[ \,(1+|\sqrt{n}\triangle
_{i}^{n}Y|^{2p}+|\beta _{i}^{n}|^{2p})\,|\,\mathcal{F}_{\frac{i-1}{n}}\right]
}\cdot \sqrt{\mathrm{E}\left[ (\sqrt{n}\triangle _{i}^{n}Y-\beta
_{i}^{n})^{4}\,|\,\mathcal{F}_{\frac{i-1}{n}}\right] } \\
&\leq &C\,\sqrt{\mathrm{E}\left[ \,\left( \int_{(i-1)/n}^{i/n}\left( \sigma
_{u-}-\sigma _{\frac{i-1}{n}}\right) \,\mathrm{d}W_{u}\right) ^{4}\,|\,%
\mathcal{F}_{\frac{i-1}{n}}\right] } \\
&\leq &C\,\sqrt{\mathrm{E}\left[ \left( \int_{(i-1)/n}^{i/n}\left( \sigma
_{u-}-\sigma _{\frac{i-1}{n}}\right) ^{2}\,\mathrm{d}u\right) ^{2}\,|\,%
\mathcal{F}_{\frac{i-1}{n}}\right] }.
\end{eqnarray*}%
Thus
\begin{eqnarray*}
\sum_{i=1}^{[nt]}\,\mathrm{E}[(\xi _{i}^{n})^{2}] &\leq &Cn\,\frac{t}{n}%
\,\sum_{i=1}^{[nt]}\mathrm{E}\,\left[ \sqrt{\mathrm{E}\,\left[ \left(
\int_{(i-1)/n}^{i/n}\left( \sigma _{u-}-\sigma _{\frac{i-1}{n}}\right) ^{2}\,%
\mathrm{d}u\right) ^{2}\,|\,\mathcal{F}_{\frac{i-1}{n}}\right] }\right] \, \\
&\leq &C\,tn\,\sqrt{\frac{1}{n}\,\sum_{i=1}^{[nt]}\mathrm{E}\left[ \left(
\int_{(i-1)/n}^{i/n}\left( \sigma _{u-}-\sigma _{\frac{i-1}{n}}\right) ^{2}\,%
\mathrm{d}u\right) ^{2}\right] } \\
&\leq &Ctn\,\sqrt{\frac{1}{n^{2}}\,\sum_{i=1}^{[nt]}\mathrm{E}\left[
\,\int_{(i-1)/n}^{i/n}\left( \sigma _{u-}-\sigma _{\frac{i-1}{n}}\right)
^{4}\,\mathrm{d}u\right] \,} \\
&\leq &Ct\,\sqrt{\,\sum_{i=1}^{[nt]}\int_{(i-1)/n}^{i/n}\mathrm{E}\left[
\left( \sigma _{u-}-\sigma _{\frac{i-1}{n}}\right) ^{2}\,\right] \mathrm{d}%
u\,} \\
&\rightarrow &\,0\vspace{2mm},
\end{eqnarray*}%
as $n\rightarrow \infty $ by Lebesgue's Theorem and the boundedness of $%
(\sigma _{t})$.

\noindent $\square $

To prove the convergence (\ref{Main result reform}) it suffices, using
Theorem \ref{theorem B} above, to prove that
\begin{equation*}
\left( U_{t}^{n}(g)-\sqrt{n}\,\left\{ \,X_{t}^{n}(g)-\int_{0}^{t}\rho
_{\sigma _{u}}(g)\,\mathrm{d}u\right\} \,\right) \overset{P}{\rightarrow }0.%
\vspace{1mm}
\end{equation*}%
But as
\begin{equation*}
U_{t}^{n}(g)-\sqrt{n}\,X_{t}^{n}(g)=-\frac{1}{\sqrt{n}}\,\sum_{i=1}^{[nt]}%
\mathrm{E}\left[ \,g(\sqrt{n}\,\triangle _{i}^{n}Y)\,|\,\mathcal{F}_{\frac{%
i-1}{n}}\right]
\end{equation*}%
and, as is easily seen,
\begin{equation*}
\left( \sqrt{n}\,\int_{0}^{t}\rho _{\sigma _{u}}\,\left( g\right) \mathrm{d}%
u-\,\sum_{i=1}^{[nt]}\sqrt{n}\,\int_{(i-1)/n}^{i/n}\rho _{\sigma _{u}}(g)\,%
\mathrm{d}u\right) \overset{P}{\rightarrow }0,\vspace{2mm}
\end{equation*}%
the job is to prove that
\begin{equation*}
\,\sum_{i=1}^{[nt]}\eta _{i}^{n}\,\overset{P}{\rightarrow }0\ \ \ \text{for
all}\ t>0,
\end{equation*}%
where for $i,n\geq 1$
\begin{equation*}
\eta _{i}^{n}=\,\frac{1}{\sqrt{n}}\,\mathrm{E}\,\left[ g(\sqrt{n}\,\triangle
_{i}^{n}Y)\,|\,\mathcal{F}_{\frac{i-1}{n}}\right] \,-\sqrt{n}%
\,\int_{(i-1)/n}^{i/n}\rho _{\sigma _{u}}(g)\,\mathrm{d}u.\vspace{1mm}
\end{equation*}%
Fix $t>0$ and write, for all $i,n\geq 1$,
\begin{equation*}
\eta _{i}^{n}=\eta (1)_{i}^{n}+\eta (2)_{i}^{n}
\end{equation*}%
where
\begin{equation}
\eta (1)_{i}^{n}=\frac{1}{\sqrt{n}}\,\,\left\{ \mathrm{E}\left[ \,g(\sqrt{n}%
\,\triangle _{i}^{n}Y)\,|\,\mathcal{F}_{\frac{i-1}{n}}\right] -\rho _{\sigma
_{\frac{i-1}{n}}}(g)\,\right\}
\end{equation}%
and
\begin{equation}
\eta (2)_{i}^{n}=\sqrt{n}\,\int_{(i-1)/n}^{i/n}\left\{ \rho _{\sigma
_{u}}(g)-\rho _{\sigma _{\frac{i-1}{n}}}(g)\right\} \mathrm{d}u.\vspace{1mm}
\end{equation}

We will now separately prove
\begin{equation}
\,\eta (1)^{n}=\sum_{i=1}^{[nt]}\eta (1)_{i}^{n}\,\overset{P}{\rightarrow }0
\label{13a}
\end{equation}%
and%
\begin{equation}
\,\eta (2)^{n}=\,\sum_{i=1}^{[nt]}\eta (2)_{i}^{n}\,\overset{P}{\rightarrow }%
0\vspace{1mm}.  \label{13b}
\end{equation}

\subsection{Some auxiliary estimates\label{subsection: intermediate limiting
results}}

In order to show (\ref{13a}) and (\ref{13b}) we need some refinements of the
estimate (\ref{8}) above. To state these we split up $(\sqrt{n}\,\triangle
_{i}^{n}Y-\beta _{i}^{n})$ into several terms. By definition
\begin{equation*}
\sqrt{n}\,\triangle _{i}^{n}Y-\beta _{i}^{n}=\sqrt{n}\,\int_{(i-1)/n}^{i/n}%
\,\left( \sigma _{u-}-\sigma _{\frac{i-1}{n}}\right) \,\mathrm{d}W_{u}%
\vspace{1mm}
\end{equation*}%
for all $i,n\geq 1$. Writing
\begin{equation*}
E_{n}=\{x\in \mathrm{E}\,|\,|\Psi (x)|>1/\sqrt{n}\,\}
\end{equation*}%
the difference $\sigma _{u}-\sigma _{\frac{i-1}{n}}$ equals
\begin{eqnarray*}
&&\int_{(i-1)/n}^{u}a_{s}^{\ast }\,ds+\int_{(i-1)/n}^{u}\sigma _{s-}^{\ast
}\,\mathrm{d}W_{s}+\int_{(i-1)/n}^{u}v_{s-}^{\ast }\,\mathrm{d}%
V_{s}+\int_{(i-1)/n}^{u}\int_{E}\phi (s-,x)\,(\mu -\nu )(\mathrm{d}s\,%
\mathrm{d}x)\vspace{1mm} \\
&=&\displaystyle\sum_{j=1}^{5}\xi (j)_{i}^{n}(u),
\end{eqnarray*}%
for $i,n\geq 1$ and $u\geq (i-1)/n$ where
\begin{eqnarray*}
\displaystyle\xi (1)_{i}^{n}(u) &=&\int_{(i-1)/n}^{u}a_{s}^{\ast }\,\mathrm{d%
}s+\int_{(i-1)/n}^{u}\,\left( \sigma _{s-}^{\ast }-\sigma _{\frac{i-1}{n}%
}^{\ast }\right) \,\mathrm{d}W_{s}+\int_{(i-1)/n}^{u}\,\left( v_{s-}^{\ast
}-v_{\frac{i-1}{n}}^{\ast }\right) \,\mathrm{d}V_{s} \\
\displaystyle\xi (2)_{i}^{n}(u) &=&\sigma _{\frac{i-1}{n}}^{\ast }\,\left(
W_{u}-W_{\frac{i-1}{n}}\right) +v_{\frac{i-1}{n}}^{\ast }\,\left( V_{u}-V_{%
\frac{i-1}{n}}\right) \\
\displaystyle\xi (3)_{i}^{n}(u) &=&\int_{(i-1)/n}^{u}\int_{E_{n}^{c}}\phi
(s-,x)\,(\mu -\nu )(\mathrm{d}s\,\mathrm{d}x) \\
\displaystyle\xi (4)_{i}^{n}(u) &=&\int_{(i-1)/n}^{u}\int_{E_{n}}\left\{
\phi (s-,x)-\phi \left( \frac{i-1}{n},x\right) \right\} \,(\mu -\nu )(%
\mathrm{d}s\,\mathrm{d}x) \\
\displaystyle\xi (5)_{i}^{n}(u) &=&\int_{(i-1)/n}^{u}\int_{E_{n}}\phi \left(
\frac{i-1}{n},x\right) \,(\mu -\nu )(\mathrm{d}s\,\mathrm{d}x)\vspace{1mm}
\end{eqnarray*}%
That is, for $i,n\geq 1$,
\begin{equation}
\sqrt{n}\,\triangle _{i}^{n}Y-\beta _{i}^{n}=\sum_{j=1}^{5}\xi (j)_{i}^{n}
\label{17}
\end{equation}%
where
\begin{equation*}
\ \xi (j)_{i}^{n}=\sqrt{n}\,\int_{(i-1)/n}^{i/n}\,\xi (j)_{i}^{n}(u-)\,%
\mathrm{d}W_{u}\ \ \ \ \mbox{\rm for}\ j=1,2,3,4,5.\vspace{1mm}
\end{equation*}%
The specific form of the variables implies, using Burkholder-Davis-Gundy
inequalities, that for every $q\geq 2$ we have
\begin{eqnarray*}
\mathrm{E}[\,|\xi (j)_{i}^{n}|^{q}\,] &\leq &C_{q}\,n^{q/2}\,\mathrm{E}\,%
\left[ \left( \int_{(i-1)/n}^{i/n}\,\xi (j)_{i}^{n}(u)^{2}\,\mathrm{d}%
u\right) ^{q/2}\right] \\
&\leq &\displaystyle n\int_{(i-1)/n}^{i/n}\,\mathrm{E}[\,|\xi
(j)_{i}^{n}(u)|^{q}\,]\mathrm{\,d}u \\
&\leq &\displaystyle\sup_{(i-1)/n\leq u\leq i/n}\,\mathrm{E}[\,|\xi
(j)_{i}^{n}(u)|^{q}\,]
\end{eqnarray*}%
for all $i,n\geq 1$ and all $j$. These terms will now be estimated. This is
done in the following series of lemmas where $i$ and $n$ are arbitrary and
we use the notation
\begin{equation*}
d_{i}^{n}=\int_{(i-1)/n}^{i/n}\mathrm{E}\,\left[ \left( \sigma _{s-}^{\ast
}-\sigma _{\frac{i-1}{n}}^{\ast }\right) ^{2}+\left( v_{s-}^{\ast }-v_{\frac{%
i-1}{n}}^{\ast }\right) ^{2}+\int_{E}\left\{ \phi (s-,x)-\phi \left( \frac{%
i-1}{n},x\right) \right\} ^{2}\,F(\mathrm{d}x)\right] \mathrm{d}s.\vspace{1mm%
}
\end{equation*}

\begin{lemma}
\label{lemma 1st}
\begin{equation*}
\mathrm{E}[\,(\xi (1)_{i}^{n})^{2}]\leq C_{1}\cdot (1/n^{2}+d_{i}^{n}).
\end{equation*}
\end{lemma}

\begin{lemma}
\label{lemma 2nd}%
\begin{equation*}
\mathrm{E}[\,(\xi (2)_{i}^{n})^{2}]\leq C_{2}/n.
\end{equation*}
\end{lemma}

\begin{lemma}
\begin{equation*}
\mathrm{E}[\,(\xi (3)_{i}^{n})^{2}]\leq C_{3}\,\varphi (1/\sqrt{n})/n,
\end{equation*}%
where%
\begin{equation*}
\varphi (\epsilon )=\int_{\{\,|\Psi |\leq \epsilon \,\}}\Psi (x)^{2}\,F(%
\mathrm{d}x).
\end{equation*}
\end{lemma}

\begin{lemma}
\begin{equation*}
\mathrm{E}[\,(\xi (4)_{i}^{n})^{2}]\leq C_{4}\,d_{i}^{n}.
\end{equation*}
\end{lemma}

\begin{lemma}
\label{lemma 5th}
\begin{equation*}
\mathrm{E}[\,(\xi (5)_{i}^{n})^{2}]\leq C_{5}/n.
\end{equation*}
\end{lemma}

The proofs of these five Lemmas rely on straightforward martingale
inequalities.

Observe that Lebesgue's Theorem ensures, since the processes involved are
assumed c\`{a}dl\`{a}g and uniformly bounded, that as $n\rightarrow \infty $
\begin{equation*}
\sum_{i=1}^{[nt]}d_{i}^{n}\,\rightarrow 0\ \ \ \ \text{for all}\ t>0.\vspace{%
1mm}
\end{equation*}

Taken together these statements imply the following result.

\begin{corollary}
\noindent \emph{For\ all\ $t>0$ as }$n\rightarrow \infty $\emph{\ }
\begin{equation*}
\sum_{i=1}^{[nt]}\,\left\{ \mathrm{E}[\,(\xi (1)_{i}^{n})^{2}]+\mathrm{E}%
[\,(\xi (3)_{i}^{n})^{2}]+\mathrm{E}[\,(\xi (4)_{i}^{n})^{2}]\right\}
\,)\,\rightarrow 0.\vspace{1mm}
\end{equation*}
\end{corollary}

Below we shall invoke this Corollary as well as Lemmas \ref{lemma 2nd} and %
\ref{lemma 5th}.\ \newline

\subsection{Proof of $\,\protect\eta (2)^{n}\protect\overset{P}{\rightarrow }%
0$ \label{subsection:13b}}

Recall we wish to show that
\begin{equation}
\,\eta (2)^{n}=\sum_{i=1}^{[nt]}\eta (2)_{i}^{n}\,\overset{P}{\rightarrow }0.
\label{eqn 7}
\end{equation}%
>From now on let $t>0$ be fixed. We split the $\eta (2)_{i}^{n}$'s according
to
\begin{equation*}
\eta (2)_{i}^{n}=\eta ^{\prime }(2)_{i}^{n}+\eta ^{\prime \prime
}(2)_{i}^{n}\ \ \ \ i,n\geq 1
\end{equation*}%
where, writing $\Phi (x)$ for $\rho _{x}(g)$,
\begin{equation*}
\eta ^{\prime }(2)_{i}^{n}=\sqrt{n}\ \Phi ^{\prime }\left( \sigma _{\frac{i-1%
}{n}}\right) \int_{(i-1)/n}^{i/n}\left( \sigma _{u}-\sigma _{\frac{i-1}{n}%
}\right) \,\mathrm{d}u
\end{equation*}%
and
\begin{equation*}
\eta ^{\prime \prime }(2)_{i}^{n}=\sqrt{n}\,\int_{(i-1)/n}^{i/n}\,\,\left\{
\Phi (\sigma _{u})-\Phi \left( \sigma _{\frac{i-1}{n}}\right) -\Phi ^{\prime
}\left( \sigma _{\frac{i-1}{n}}\right) \cdot \left( \sigma _{u}-\sigma _{%
\frac{i-1}{n}}\right) \right\} \,\,\mathrm{d}u.
\end{equation*}%
Observe that the assumptions on $g$ imply that $x\mapsto \Phi (x)$ is
differentiable with a bounded derivative on any bounded interval not
including $0$; in particular\thinspace (see (H2a))
\begin{equation}
|\,\Phi (x)-\Phi (y)-\Phi ^{\prime }(y)\cdot (x-y)\,|\leq \Psi (|x-y|)\cdot
|x-y|,\ \ \ x^{2},y^{2}\in (a,b),  \label{eqn 8}
\end{equation}%
where $\Psi :\mathbf{R}_{+}\rightarrow \mathbf{R}_{+}$ is continuous,
increasing and $\Psi (0)=0$. \vspace{1mm}

With this notation we shall prove (\ref{eqn 7}) by showing%
\begin{equation*}
\,\sum_{i=1}^{[nt]}\eta ^{\prime }(2)_{i}^{n}\,\overset{P}{\rightarrow }0
\end{equation*}%
and
\begin{equation*}
\ \,\sum_{i=1}^{[nt]}\eta ^{\prime \prime }(2)_{i}^{n}\,\overset{P}{%
\rightarrow }0.\vspace{1mm}
\end{equation*}%
Inserting the description of $(\sigma _{t})$\thinspace (see (H1a)) we may
write
\begin{equation*}
\eta ^{\prime }(2)_{i}^{n}=\eta ^{\prime }(2,1)_{i}^{n}+\eta ^{\prime
}(2,2)_{i}^{n}
\end{equation*}%
where for all $i,n\geq 1$
\begin{equation*}
\eta ^{\prime }(2,1)_{i}^{n}=\sqrt{n}\ \Phi ^{\prime }\left( \sigma _{\frac{%
i-1}{n}}\right) \int_{(i-1)/n}^{i/n}\left( \int_{(i-1)/n}^{u}a_{s}^{\ast }\,%
\mathrm{d}s\right) \,\,\mathrm{d}u
\end{equation*}%
and
\begin{eqnarray*}
\eta ^{\prime }(2,2)_{i}^{n} &=&\displaystyle\sqrt{n}\ \Phi ^{\prime }\left(
\sigma _{\frac{i-1}{n}}\right) \int_{(i-1)/n}^{i/n}\,\left[
\int_{(i-1)/n}^{u}\,\sigma _{s-}^{\ast }\,\mathrm{d}W_{s}+\int_{(i-1)/n}^{u}%
\,v_{s-}^{\ast }\,\mathrm{d}V_{s}\right. \, \\
&&+\displaystyle\left. \int_{E}\phi (s-,x)\,(\mu -\nu )(\mathrm{d}s\,\mathrm{%
d}x)\right] \mathrm{d}u.\vspace{1mm}
\end{eqnarray*}%
By (H2a) and (\ref{eqn 8}) and the uniform boundedness of $(a_{t}^{\ast })$
we have
\begin{equation*}
|\eta ^{\prime }(2,1)_{i}^{n}|\leq C\,\sqrt{n}\,\int_{(i-1)/n}^{i/n}\left\{
u-(i-1)/n\right\} \,\mathrm{d}u\leq C/n^{3/2}
\end{equation*}%
for all $i,n\geq 1$ and thus
\begin{equation*}
\,\sum_{i=1}^{[nt]}\eta ^{\prime }(2,1)_{i}^{n}\,\overset{P}{\rightarrow }0.%
\vspace{1mm}
\end{equation*}%
Since
\begin{equation*}
(W_{t}),\ (V_{t})\ \text{and}\ \left( \int_{0}^{t}\int_{E}\phi (s-,x)(\mu
-\nu )(\mathrm{d}s\,\mathrm{d}x)\right) \vspace{1mm}
\end{equation*}%
are all martingales we have
\begin{equation*}
\mathrm{E}\left[ \eta ^{\prime }(2,2)_{i}^{n}\,|\,\mathcal{F}_{\frac{i-1}{n}}%
\right] =0\ \ \ \text{for all}\quad i,n\geq 1.\vspace{1mm}
\end{equation*}

By Doob's inequality it is therefore feasible to estimate
\begin{equation*}
\sum_{i=1}^{[nt]}\,\mathrm{E}[\,(\eta ^{\prime }(2,2)_{i}^{n})^{2}].\vspace{%
1mm}
\end{equation*}%
Inserting again the description of $(\sigma _{t})$ we find, applying simple
inequalities, in particular Jensen's, that
\begin{eqnarray*}
&&(\eta ^{\prime }(2,2)_{i}^{n})^{2} \\
&\leq &\displaystyle C\,n\,\left( \int_{(i-1)/n}^{i/n}\left\{
\int_{(i-1)/n}^{u}\,\sigma _{s-}^{\ast }\,\mathrm{d}W_{s}\right\} \,\mathrm{d%
}u\right) ^{2}+C\,n\,\left( \,\int_{(i-1)/n}^{i/n}\left\{
\int_{(i-1)/n}^{u}\,v_{s-}^{\ast }\,\mathrm{d}V_{s}\right\} \,\mathrm{d}%
u\right) ^{2} \\
&&\displaystyle+C\,n\,\left( \,\int_{(i-1)/n}^{i/n}\int_{(i-1)/n}^{u}\left\{
\int_{E}\phi (s-,x)\,(\mu -\nu )(\mathrm{d}s\,\mathrm{d}x)\right\} \mathrm{d}%
u\right) ^{2} \\
&\leq &\displaystyle C\,\int_{(i-1)/n}^{i/n}\left(
\,\int_{(i-1)/n}^{u}\,\sigma _{s-}^{\ast }\,\mathrm{d}W_{s}\,\right) ^{2}\,%
\mathrm{d}u+C\,\int_{(i-1)/n}^{i/n}\left( \,\int_{(i-1)/n}^{u}\,v_{s-}^{\ast
}\,\mathrm{d}V_{s}\,\right) ^{2}\,\mathrm{d}u \\
&&\displaystyle+C\,\int_{(i-1)/n}^{i/n}\left(
\,\int_{(i-1)/n}^{u}\int_{E}\phi (s-,x)\,(\mu -\nu )(\mathrm{d}s\,\mathrm{d}%
x)\,\right) ^{2}\,\mathrm{d}u.\vspace{1mm}
\end{eqnarray*}%
The properties of the Wiener integrals and the uniform boundedness of $%
(\sigma _{t}^{\ast })$ and $(v_{t}^{\ast })$ ensure that
\begin{equation*}
\mathrm{E}\left[ \left( \,\int_{(i-1)/n}^{u}\,\sigma _{s-}^{\ast }\,\mathrm{d%
}W_{s}\,\right) ^{2}\,|\,\mathcal{F}_{\frac{i-1}{n}}\right] \leq C\cdot
\left( u-\frac{i-1}{n}\right)
\end{equation*}%
and likewise
\begin{equation*}
\mathrm{E}\left[ \left( \,\int_{(i-1)/n}^{u}\,v_{s-}^{\ast }\,\mathrm{d}%
V_{s}\,\right) ^{2}\,|\,\mathcal{F}_{\frac{i-1}{n}}\right] \leq C\cdot
\left( u-\frac{i-1}{n}\right) \vspace{1mm}
\end{equation*}%
for all $i,n\geq 1$. Likewise for the Poisson part we have
\begin{eqnarray*}
&&\displaystyle\mathrm{E}\left[ \left( \,\int_{(i-1)/n}^{u}\int_{E}\phi
(s-,x)\,(\mu -\nu )(\mathrm{d}s\,\mathrm{d}x)\,\right) ^{2}\,|\,\mathcal{F}_{%
\frac{i-1}{n}}\right] \\
&\leq &\displaystyle C\int_{(i-1)/n}^{u}\int_{E}\mathrm{E}[\phi
^{2}(s,x)\,|\,\mathcal{F}_{\frac{i-1}{n}}]\,F(\mathrm{d}x)\,\mathrm{d}s%
\vspace{1mm}
\end{eqnarray*}%
yielding a similar bound. Putting all this together we have for all $i,n\geq
1$
\begin{eqnarray*}
\mathrm{E}[\,(\eta ^{\prime }(2,2)_{i}^{n})^{2}\,|\,\mathcal{F}_{\frac{i-1}{n%
}}] &\leq &C\,\int_{(i-1)/n}^{i/n}(u-(i-1)/n)\,\mathrm{d}u \\
&\leq &C/n^{2}.
\end{eqnarray*}%
Thus as $n\rightarrow \infty $ so
\begin{equation*}
\sum_{i=1}^{[nt]}\mathrm{E}[\,(\eta ^{\prime }(2,2)_{i}^{n})^{2}]\rightarrow
0.\vspace{1mm}
\end{equation*}%
and since
\begin{equation*}
\mathrm{E}\left[ \eta ^{\prime }(2,2)_{i}^{n}\,|\,\mathcal{F}_{\frac{i-1}{n}}%
\right] =0\ \ \ \ \text{for all}\quad i,n\geq 1\vspace{1mm}
\end{equation*}%
we deduce from Doob's inequality that
\begin{equation*}
\,\sum_{i=1}^{[nt]}\eta ^{\prime }(2,2)_{i}^{n}\,\overset{P}{\rightarrow }0
\end{equation*}%
proving\ altogether
\begin{equation*}
\,\sum_{i=1}^{[nt]}\eta ^{\prime }(2)_{i}^{n}\,\overset{P}{\rightarrow }0.%
\vspace{1mm}
\end{equation*}%
Applying once more (H2a) and (\ref{eqn 8}) we have for every $\epsilon >0$
and every $i,n$ that
\begin{eqnarray*}
|\eta ^{\prime \prime }(2)_{i}^{n}| &\leq &\sqrt{n}\int_{(i-1)/n}^{i/n}\,%
\Psi \left( \left\vert \sigma _{u}-\sigma _{\frac{i-1}{n}}\right\vert
\right) \cdot \left\vert \sigma _{u}-\sigma _{\frac{i-1}{n}}\right\vert \,%
\mathrm{d}u \\
&\leq &\sqrt{n}\,\Psi (\epsilon )\int_{(i-1)/n}^{i/n}\,\left\vert \sigma
_{u}-\sigma _{\frac{i-1}{n}}\right\vert \,\mathrm{d}u+\sqrt{n}\,\Psi (2\sqrt{%
b})/\epsilon \int_{(i-1)/n}^{i/n}\,\left\vert \sigma _{u}-\sigma _{\frac{i-1%
}{n}}\right\vert ^{2}\,\mathrm{d}u.\vspace{1mm}
\end{eqnarray*}%
Thus from (\ref{8}) and its consequence
\begin{equation*}
\mathrm{E}\,\left[ \left\vert \sigma _{u}-\sigma _{\frac{i-1}{n}}\right\vert
\,\right] \leq C/\sqrt{n}
\end{equation*}%
we get
\begin{equation*}
\sum_{i=1}^{[nt]}\mathrm{E}[\,|\eta ^{\prime \prime }(2)_{i}^{n}|\,]\leq
Ct\,\Psi (\epsilon )+\frac{C\,\Psi (b)}{\sqrt{n}\,\epsilon }
\end{equation*}%
for all $n$ and all $\epsilon $. Letting here first $n\rightarrow \infty $
and then $\epsilon \rightarrow 0$ we may conclude that as $n\rightarrow
\infty $
\begin{equation*}
\sum_{i=1}^{[nt]}\mathrm{E}[\,|\eta ^{\prime \prime
}(2)_{i}^{n}|\,]\rightarrow 0\
\end{equation*}%
implying the convergence
\begin{equation*}
\,\sum_{i=1}^{[nt]}\eta (2)_{i}^{n}\,\overset{P}{\rightarrow }0.\vspace{1mm}
\end{equation*}%
Thus ending the proof of (\ref{13b}).

\noindent $\square $

\subsection{Proof of $\protect\eta (1)^{n}\protect\overset{P}{\rightarrow }0$
\label{subsection:proof of 13a}}

Recall we are to show that \textbf{\ }
\begin{equation}
\eta (1)^{n}=\,\sum_{i=1}^{[nt]}\eta (1)_{i}^{n}\,\overset{P}{\rightarrow }0.
\end{equation}%
Let still $t>0$ be fixed. Recall that
\begin{eqnarray*}
\eta (1)_{i}^{n} &=&\frac{1}{\sqrt{n}}\,\left\{ \mathrm{E}\left[ \,g(\sqrt{n}%
\,\triangle _{i}^{n}Y)\,|\,\mathcal{F}_{\frac{i-1}{n}}\right] \,-\rho
_{\sigma _{\frac{i-1}{n}}}(g)\right\} \\
&=&\frac{1}{\sqrt{n}}\,\mathrm{E}\,\left[ g(\sqrt{n}\,\triangle
_{i}^{n}Y)-g(\beta _{i}^{n})\,|\,\mathcal{F}_{\frac{i-1}{n}}\right] .\vspace{%
1mm}
\end{eqnarray*}%
Introduce the notation\thinspace (recall the assumption (K))
\begin{equation*}
A_{i}^{n}=\{\,|\sqrt{n}\,\triangle _{i}^{n}Y-\beta _{i}^{n}|>\,d(\beta
_{i}^{n},B)/2\,\}.\vspace{1mm}
\end{equation*}%
Since $B$ is a Lebesgue null set and $\beta _{i}^{n}$ is absolutely
continuous, $g^{\prime }(\beta _{i}^{n})$ is defined $a.s.$ and, by
assumption, $g$ is differentiable on the interval joining $\triangle
_{i}^{n}Y(\omega )$ and $\beta _{i}^{n}(\omega )$ for all $\omega \in
A_{i}^{n\,c}$. Thus, using the Mean Value Theorem, we may for all $i,n\geq 1$
write
\begin{eqnarray*}
&&g(\sqrt{n}\,\triangle _{i}^{n}Y)-g(\beta _{i}^{n}) \\
&=&\left\{ g(\sqrt{n}\,\triangle _{i}^{n}Y)-g(\beta _{i}^{n})\right\} \cdot
\mathbf{1}_{A_{i}^{n}} \\
&&+g^{\prime }(\beta _{i}^{n})\cdot (\sqrt{n}\,\triangle _{i}^{n}Y-\beta
_{i}^{n})\cdot \mathbf{1}_{A_{i}^{n\,c}} \\
&&+\left\{ g^{\prime }(\alpha _{i}^{n})-g^{\prime }(\beta _{i}^{n})\right\}
\cdot (\sqrt{n}\,\triangle _{i}^{n}Y-\beta _{i}^{n})\cdot \mathbf{1}%
_{A_{i}^{n\,c}} \\
&=&\sqrt{n}\,\left\{ \delta (1)_{i}^{n}+\delta (2)_{i}^{n}+\delta
(3)_{i}^{n}\right\} ,\vspace{1mm}
\end{eqnarray*}%
where $\alpha _{i}^{n}$ are random points lying in between $\sqrt{n}%
\,\triangle _{i}^{n}Y$ and $\beta _{i}^{n}$, i.e.
\begin{equation*}
\sqrt{n}\,\triangle _{i}^{n}Y\wedge \beta _{i}^{n}\leq \alpha _{i}^{n}\leq
\sqrt{n}\,\triangle _{i}^{n}Y\vee \beta _{i}^{n},
\end{equation*}%
and%
\begin{equation*}
\begin{array}{lll}
\delta (1)_{i}^{n} & = & \left[ \,\left\{ g(\sqrt{n}\,\triangle
_{i}^{n}Y)-g(\beta _{i}^{n})\right\} -g^{\prime }(\beta _{i}^{n})\cdot (%
\sqrt{n}\,\triangle _{i}^{n}Y-\beta _{i}^{n})\,\right] \cdot \mathbf{1}%
_{A_{i}^{n}}/\sqrt{n} \\
\delta (2)_{i}^{n} & = & \left\{ g^{\prime }(\alpha _{i}^{n})-g^{\prime
}(\beta _{i}^{n})\right\} \cdot (\sqrt{n}\,\triangle _{i}^{n}Y-\beta
_{i}^{n})\cdot \mathbf{1}_{A_{i}^{n\,c}}/\sqrt{n} \\
\delta (3)_{i}^{n} & = & g^{\prime }(\beta _{i}^{n})\cdot (\sqrt{n}%
\,\triangle _{i}^{n}Y-\beta _{i}^{n})/\sqrt{n}.%
\end{array}%
\end{equation*}%
Thus it suffices to prove
\begin{equation*}
\,\sum_{i=1}^{[nt]}\mathrm{E}\,\left[ \delta (k)_{i}^{n}\,|\,\mathcal{F}_{%
\frac{i-1}{n}}\right] \,\,\overset{P}{\rightarrow }0,\ \ \ k=1,2,3.
\end{equation*}

Consider the case $k=1$. Using (K) and the fact that $\beta _{i}^{n}$ is
absolutely continuous we have a.s.%
\begin{eqnarray*}
&&|g(\sqrt{n}\,\triangle _{i}^{n}Y)-g(\beta _{i}^{n})| \\
&\leq &M(1+|\sqrt{n}\,\triangle _{i}^{n}Y-\beta _{i}^{n}|^{p}+|\beta
_{i}^{n}|^{p})\cdot |\sqrt{n}\,\triangle _{i}^{n}Y-\beta _{i}^{n}| \\
&\leq &(2^{p}+1)M(1+|\sqrt{n}\,\triangle _{i}^{n}Y|^{p}+|\beta
_{i}^{n}|^{p})\cdot |\sqrt{n}\,\triangle _{i}^{n}Y-\beta _{i}^{n}|,
\end{eqnarray*}%
and
\begin{equation*}
|\,g^{\prime }(\beta _{i}^{n})\cdot (\sqrt{n}\,\triangle _{i}^{n}Y-\beta
_{i}^{n})\,|\leq M(1+|\beta _{i}^{n}|^{p})\cdot |\sqrt{n}\,\triangle
_{i}^{n}Y-\beta _{i}^{n}|.\vspace{1mm}
\end{equation*}%
By Cauchy-Schwarz's inequality $\mathrm{E}[\,|\delta (1)_{i}^{n}|\,]$ is
therefore for all $i,n\geq 1$ less than
\begin{equation*}
C\cdot \mathrm{E}[\,1+|\sqrt{n}\,\triangle _{i}^{n}Y|^{3p}+|\beta
_{i}^{n}|^{3p}]^{1/3}\cdot \mathrm{E}[\,(\sqrt{n}\,\triangle _{i}^{n}Y-\beta
_{i}^{n})^{2}/n\,]^{1/2}\cdot P(A_{i}^{n})^{1/6}\vspace{1mm}
\end{equation*}%
implying for fixed $t,$ by means of (\ref{2}), that
\begin{eqnarray*}
\mathrm{E}[\left[ \sum_{i=1}^{[nt]}|\,\delta (1)_{i}^{n}|\right] \, &\leq
&C\cdot \sup_{i\geq 1}P(A_{i}^{n})^{1/6}\,\sum_{i=1}^{[nt]}\mathrm{E}%
[\,(\triangle _{i}^{n}Y-\beta _{i}^{n})^{2}/n\,]^{1/2}\vspace{1mm} \\
&\leq &C\cdot \sup_{i\geq 1}P(A_{i}^{n})^{1/6}\,\sum_{i=1}^{[nt]}1/n \\
&\leq &Ct\cdot \sup_{i\geq 1}P(A_{i}^{n})^{1/6}.\vspace{1mm}
\end{eqnarray*}%
For all $i,n\geq 1$ we have for every $\epsilon >0$
\begin{eqnarray*}
P(A_{i}^{n}) &\leq &P(A_{i}^{n}\cap \{d(\beta _{i}^{n},B)\leq \epsilon
\})+P(A_{i}^{n}\cap \{d(\beta _{i}^{n},B)>\epsilon \})\vspace{1mm} \\
&\leq &P(d(\beta _{i}^{n},B)\leq \epsilon )+P(|\sqrt{n}\,\triangle
_{i}^{n}Y-\beta _{i}^{n}|>\epsilon /2)\vspace{1mm} \\
&\leq &P(d(\beta _{i}^{n},B)\leq \epsilon )+\frac{4}{\epsilon ^{2}}\cdot
\mathrm{E}[\,(\sqrt{n}\,\triangle _{i}^{n}Y-\beta _{i}^{n})^{2}]\vspace{1mm}
\\
&\leq &P(d(\beta _{i}^{n},B)\leq \epsilon )+\frac{C}{n\,\epsilon ^{2}}.
\end{eqnarray*}%
But (H2a) implies that the densities of $\beta _{i}^{n}$ are pointwise
dominated by a Lebesgue integrable function $h_{a,b}$ providing, for all $%
i,n\geq 1$, the estimate
\begin{eqnarray}
P(A_{i}^{n}) &\leq &\int_{\{x\,|\,d(x,B)\leq \epsilon \}}h_{a,b}\,\mathrm{d}%
\lambda _{1}+\frac{C}{n\,\epsilon ^{2}}  \label{eqn 10} \\
&=&\alpha _{\epsilon }+\frac{C}{n\,\epsilon ^{2}}.\vspace{1mm}  \notag
\end{eqnarray}%
Observe $\lim_{\epsilon \rightarrow 0}\alpha _{\epsilon }=0$. Taking now in (%
\ref{eqn 10}) $\sup $ over $i$ and then letting first $n\rightarrow \infty $
and then $\epsilon \downarrow 0$ we get
\begin{equation*}
\lim_{n}\,\sup_{i\geq 1}\,P(A_{i}^{n})=0
\end{equation*}%
proving that
\begin{equation*}
\mathrm{E}\left[ \,\sum_{i=1}^{[nt]}|\,\delta (1)_{i}^{n}|\right]
\,\rightarrow 0
\end{equation*}%
and\ thus
\begin{equation*}
\,\sum_{i=1}^{[nt]}\mathrm{E}\left[ \,\delta (1)_{i}^{n}\,|\,\mathcal{F}_{%
\frac{i-1}{n}}\right] \,\overset{P}{\rightarrow }0.\vspace{1mm}
\end{equation*}

Consider next the case $k=2$. As assumed in (K), $g$ is continuously
differentiable outside of $B$. Thus for each $A>1$ and $\epsilon >0$ there
exists a function $G_{A,\,\epsilon }:(0,1)\rightarrow \mathbf{R}_{+}$ such
that for given $0<\epsilon ^{\prime }<\epsilon /2$
\begin{equation*}
\left\vert g^{\prime }(x+y)-g^{\prime }(x)\right\vert \leq G_{A,\,\epsilon
}(\epsilon ^{\prime })\ \ \text{for all}\ |x|\leq A,\ |y|\leq \epsilon
^{\prime }<\epsilon <d(x,B).\vspace{1mm}
\end{equation*}%
Observe that $\lim_{\epsilon ^{\prime }\downarrow 0}G_{A,\,\epsilon
}(\epsilon ^{\prime })=0$ for all $A$ and $\epsilon $.\vspace{1mm} Fix $A>1$
and $\epsilon \in (0,1)$. For all $i,n\geq 1$ we have
\begin{eqnarray*}
&&|g^{\prime }(\alpha _{i}^{n})-g^{\prime }(\beta _{i}^{n})|\cdot \mathbf{1}%
_{A_{i}^{n\,c}} \\
&=&\displaystyle|g^{\prime }(\alpha _{i}^{n})-g^{\prime }(\beta
_{i}^{n})|\cdot \mathbf{1}_{A_{i}^{n\,c}}\,(\mathbf{1}_{\{|\alpha
_{i}^{n}|+|\beta _{i}^{n}|>A\}}+\mathbf{1}_{\{|\alpha _{i}^{n}|+|\beta
_{i}^{n}|\leq A\}})\vspace{1mm} \\
&\leq &\displaystyle|g^{\prime }(\alpha _{i}^{n})-g^{\prime }(\beta
_{i}^{n})|\cdot \frac{|\alpha _{i}^{n}|+|\beta _{i}^{n}|}{A}+|g^{\prime
}(\alpha _{i}^{n})-g^{\prime }(\beta _{i}^{n})|\cdot \mathbf{1}%
_{A_{i}^{n\,c}\,\cap \,\{|\alpha _{i}^{n}|+|\beta _{i}^{n}|\leq A\}} \\
&\leq &\displaystyle\frac{C}{A}\cdot (1+|\alpha _{i}^{n}|^{p}+|\beta
_{i}^{n}|^{p})^{2}+|g^{\prime }(\alpha _{i}^{n})-g^{\prime }(\beta
_{i}^{n})|\cdot \mathbf{1}_{A_{i}^{n\,c}\,\cap \,\{|\alpha _{i}^{n}|+|\beta
_{i}^{n}|\leq A\}} \\
&\leq &\displaystyle\frac{C}{A}\cdot (1+|\sqrt{n}\,\triangle
_{i}^{n}Y|^{2p}+|\beta _{i}^{n}|^{2p})+|g^{\prime }(\alpha
_{i}^{n})-g^{\prime }(\beta _{i}^{n})|\cdot \mathbf{1}_{A_{i}^{n\,c}\,\cap
\,\{|\alpha _{i}^{n}|+|\beta _{i}^{n}|\leq A\}}.\vspace{1mm}
\end{eqnarray*}%
Now writing
\begin{eqnarray*}
1 &=&\mathbf{1}_{\{d(\beta _{i}^{n},B)\leq \epsilon \}}+\mathbf{1}%
_{\{d(\beta _{i}^{n},B)>\epsilon \}}\vspace{1mm} \\
&=&\mathbf{1}_{\{d(\beta _{i}^{n},B)\leq \epsilon \}} \\
&&+\mathbf{1}_{\{d(\beta _{i}^{n},B)>\epsilon \}\,\cap \,\{|\alpha
_{i}^{n}-\beta _{i}^{n}|\leq \epsilon ^{\prime }\}} \\
&&+\mathbf{1}_{\{d(\beta _{i}^{n},B)>\epsilon \}\,\cap \,\{|\alpha
_{i}^{n}-\beta _{i}^{n}|>\epsilon ^{\prime }\}}
\end{eqnarray*}%
for all $0<\epsilon ^{\prime }<\epsilon /2$ we have
\begin{eqnarray*}
\mathbf{1}_{A_{i}^{n\,c}\,\cap \,\{|\alpha _{i}^{n}|+|\beta _{i}^{n}|\leq
A\}} &\leq &\mathbf{1}_{\{d(\beta _{i}^{n},B)\leq \epsilon \}\,\cap
\,A_{i}^{n\,c}\,\cap \,\{|\alpha _{i}^{n}|+|\beta _{i}^{n}|\leq A\}} \\
&&+\mathbf{1}_{A_{i}^{n\,c}\,\cap \,\{|\alpha _{i}^{n}|+|\beta _{i}^{n}|\leq
A\}\,\cap \,\{d(\beta _{i}^{n},B)>\epsilon \}\,\cap \,\{|\alpha
_{i}^{n}-\beta _{i}^{n}|\leq \epsilon ^{\prime }\}} \\
&&+\mathbf{1}_{A_{i}^{n\,c}\,\cap \,\{|\alpha _{i}^{n}|+|\beta _{i}^{n}|\leq
A\}\,\cap \,\{d(\beta _{i}^{n},B)>\epsilon \}}\cdot \frac{|\alpha
_{i}^{n}-\beta _{i}^{n}|}{\epsilon ^{\prime }}.
\end{eqnarray*}%
Combining this with the fact that
\begin{eqnarray*}
|g^{\prime }(\alpha _{i}^{n})-g^{\prime }(\beta _{i}^{n})| &\leq
&C(1+|\alpha _{i}^{n}|^{p}+|\beta _{i}^{n}|^{p}) \\
&\leq &CA^{p}
\end{eqnarray*}%
on $A_{i}^{n\,c}\,\cap \,\{|\alpha _{i}^{n}|+|\beta _{i}^{n}|\leq A\}$ we
obtain that
\begin{eqnarray*}
&&|g^{\prime }(\alpha _{i}^{n})-g^{\prime }(\beta _{i}^{n})|\cdot \mathbf{1}%
_{A_{i}^{n\,c}\,\cap \,\{|\alpha _{i}^{n}|+|\beta _{i}^{n}|\leq A\}} \\
&\leq &CA^{p}\cdot \left( \,\mathbf{1}_{\{d(\beta _{i}^{n},B)\leq \epsilon
\}}+\frac{|\alpha _{i}^{n}-\beta _{i}^{n}|}{\epsilon ^{\prime }}\right)
\,+G_{A,\,\epsilon }(\epsilon ^{\prime })\vspace{1mm} \\
&\leq &CA^{p}\cdot (\,\mathbf{1}_{\{d(\beta _{i}^{n},B)\leq \epsilon \}}+%
\frac{|\sqrt{n}\,\triangle _{i}^{n}Y-\beta _{i}^{n}|}{\epsilon ^{\prime }}%
\,)+G_{A,\,\epsilon }(\epsilon ^{\prime }).\vspace{1mm}
\end{eqnarray*}

Putting this together means that
\begin{eqnarray*}
\sqrt{n}\,|\delta (2)_{i}^{n}| &=&|g^{\prime }(\alpha _{i}^{n})-g^{\prime
}(\beta _{i}^{n})|\cdot |\sqrt{n}\,\triangle _{i}^{n}Y-\beta _{i}^{n}|\cdot
\mathbf{1}_{A_{i}^{n\,c}} \\
&\leq &\left\{ \frac{C}{A}\cdot (1+|\sqrt{n}\,\triangle
_{i}^{n}Y|^{2p}+|\beta _{i}^{n}|^{2p})+G_{A,\,\epsilon }(\epsilon ^{\prime
})\right\} \cdot |\sqrt{n}\,\triangle _{i}^{n}Y-\beta _{i}^{n}|\vspace{1mm}
\\
&&+\,CA^{p}\cdot \left( \mathbf{1}_{\{d(\beta _{i}^{n},B)\leq \epsilon
\}}\cdot |\sqrt{n}\,\triangle _{i}^{n}Y-\beta _{i}^{n}|+\frac{|\sqrt{n}%
\,\triangle _{i}^{n}Y-\beta _{i}^{n}|^{2}}{\epsilon ^{\prime }}\right) .%
\vspace{1mm}
\end{eqnarray*}%
Exploiting here the inequalities (\ref{2}) and (\ref{3}) we obtain, for all $%
A>1$ and $0<2\epsilon ^{\prime }<\epsilon <1$ and all $i,n\geq 1$, using H%
\"{o}lder's inequality, the following estimate
\begin{equation*}
\mathrm{E}[\,|\delta (2)_{i}^{n}|\,]\leq C\left( \frac{1}{A\,n}+\frac{%
G_{A,\,\epsilon }(\epsilon ^{\prime })}{n}+\frac{A^{p}\,\sqrt{\alpha
_{\epsilon }}}{n}+\frac{A^{p}}{\epsilon ^{\prime }\,n^{3/2}}\right) \vspace{%
1mm}
\end{equation*}%
implying for all $n\geq 1$ and $t\geq 0$ that
\begin{equation*}
\sum_{i=1}^{[nt]}\mathrm{E}[\,|\delta (2)_{i}^{n}|\,]\leq Ct\left( \frac{1}{A%
}+G_{A,\,\epsilon }(\epsilon ^{\prime })+A^{p}\,\sqrt{\alpha _{\epsilon }}+%
\frac{A^{p}}{\epsilon ^{\prime }\,n^{1/2}}\right) .\vspace{1mm}
\end{equation*}%
Choosing in this estimate first $A$ sufficiently big, then $\epsilon $
small\thinspace (recall that $\lim_{\epsilon \rightarrow 0}\alpha _{\epsilon
}=0$\thinspace ) and finally $\epsilon ^{\prime }$ small, exploiting that $%
\lim_{\epsilon ^{\prime }\downarrow 0}G_{A,\,\epsilon }(\epsilon ^{\prime
})=0$ for all $A$ and $\epsilon $, we may conclude that
\begin{equation*}
\lim_{n}\,\sum_{i=1}^{[nt]}\mathrm{E}\left[ \,|\delta (2)_{i}^{n}|\,\right]
=0
\end{equation*}%
and thus
\begin{equation*}
\sum_{i=1}^{[nt]}\mathrm{E}\left[ \,\delta (2)_{i}^{n}\,|\,\mathcal{F}_{%
\frac{i-1}{n}}\right] \,\overset{P}{\rightarrow }0.
\end{equation*}

So what remains to be proved is the convergence
\begin{equation*}
\,\sum_{i=1}^{[nt]}\mathrm{E}\,\left[ \delta (3)_{i}^{n}\,|\,\mathcal{F}_{%
\frac{i-1}{n}}\right] \,\,\overset{P}{\rightarrow }0.
\end{equation*}%
As introduced in (\ref{17})
\begin{equation*}
\sqrt{n}\,\triangle _{i}^{n}Y-\beta _{i}^{n}=\sum_{j=1}^{5}\xi
(j)_{i}^{n}=\psi (1)_{i}^{n}+\psi (2)_{i}^{n}
\end{equation*}%
for all $i,n\geq 1$ where
\begin{equation*}
\psi (1)_{i}^{n}=\xi (1)_{i}^{n}+\xi (3)_{i}^{n}+\xi (4)_{i}^{n},
\end{equation*}%
\begin{equation*}
\psi (2)_{i}^{n}=\xi (2)_{i}^{n}+\xi (5)_{i}^{n},
\end{equation*}%
and as
\begin{equation*}
\delta (3)_{i}^{n}=g^{\prime }(\beta _{i}^{n})\cdot (\psi (1)_{i}^{n}+\psi
(2)_{i}^{n})/\sqrt{n}\vspace{1mm}
\end{equation*}%
it suffices to prove
\begin{equation*}
\left( \,\sum_{i=1}^{[nt]}\mathrm{E}\left[ \,g^{\prime }(\beta
_{i}^{n})\cdot \psi (k)_{i}^{n}\,|\,\mathcal{F}_{\frac{i-1}{n}}\right] \,\,/%
\sqrt{n}\,\right) \overset{P}{\rightarrow }0,\ \ \ k=1,2.\vspace{1mm}
\end{equation*}

The case $k=1$ is handled by proving
\begin{equation}
\frac{1}{\sqrt{n}}\,\sum_{i=1}^{[nt]}\mathrm{E}[\,|g^{\prime }(\beta
_{i}^{n})\cdot \xi (j)_{i}^{n}|\,]\rightarrow 0,\ \ \ j=1,3,4.\vspace{1mm}
\label{eqn 11}
\end{equation}%
Using Jensen's inequality it is easily seen that for $j=1,3,4$
\begin{equation*}
\frac{1}{\sqrt{n}}\,\sum_{i=1}^{[nt]}\mathrm{E}[\,|g^{\prime }(\beta
_{i}^{n})\cdot \xi (j)_{i}^{n}|\,]\leq C\,t\cdot \sqrt{\frac{1}{n}%
\,\sum_{i=1}^{[nt]}\mathrm{E}[\,g^{\prime }(\beta _{i}^{n})^{2}]}\,\cdot \,%
\sqrt{\sum_{i=1}^{[nt]}\mathrm{E}[\,(\xi (j)_{i}^{n})^{2}]}
\end{equation*}%
and so using (\ref{12})
\begin{equation*}
\frac{1}{\sqrt{n}}\,\sum_{i=1}^{[nt]}\mathrm{E}[\,|g^{\prime }(\beta
_{i}^{n})\cdot \xi (j)_{i}^{n}|\,]\leq C\,t\cdot \,\sqrt{\sum_{i=1}^{[nt]}%
\mathrm{E}[\,(\xi (j)_{i}^{n})^{2}]}
\end{equation*}%
since almost surely
\begin{equation*}
|g^{\prime }(\beta _{i}^{n})|\leq C\,(1+|\beta _{i}^{n}|^{p})
\end{equation*}%
for all $i,n\geq 1$. From here, (\ref{eqn 11}) is an immediate consequence
of Lemmas \ref{lemma 1st}-\ref{lemma 5th}.\vspace{1mm}

The remaining case $k=2$ is different. The definition of $\psi (2)_{i}^{n}$
implies, using basic stochastic calculus, that $\psi (2)_{i}^{n}/\sqrt{n}$,
for all $i,n\geq 1$, may be written as
\begin{eqnarray*}
&&\int_{(i-1)/n}^{i/n}\left\{ \sigma _{\frac{i-1}{n}}^{\prime }\,\left(
W_{u}-W_{\frac{i-1}{n}}\right) +M(n,i)_{u}\right\} \,\mathrm{d}W_{u} \\
&=&\sigma _{\frac{i-1}{n}}^{\prime }\,\int_{(i-1)/n}^{i/n}\left( W_{u}-W_{%
\frac{i-1}{n}}\right) \,\mathrm{d}W_{u} \\
&&+\triangle _{i}^{n}M(n,i)\cdot \triangle _{i}^{n}W \\
&&+\int_{(i-1)/n}^{i/n}\left( W_{u}-W_{\frac{i-1}{n}}\right) \,\mathrm{d}%
M(n,i)_{u},
\end{eqnarray*}%
where $(M(n,i)_{t})$ is the martingale defined by $M(n,i)_{t}\equiv 0$ for $%
t\leq (i-1)/n$ and
\begin{equation*}
M(n,i)_{t}=v_{\frac{i-1}{n}}^{\ast }\,\left( V_{t}-V_{\frac{i-1}{n}}\right)
+\int_{(i-1)/n}^{t}\int_{E_{n}}\phi \left( \frac{i-1}{n},x\right) (\mu -\nu
)(\mathrm{d}s\,\mathrm{d}x)\vspace{1mm}
\end{equation*}%
otherwise. Thus for fixed $i,n\geq 1$
\begin{equation*}
\mathrm{E}\left[ \,g^{\prime }(\beta _{i}^{n})\cdot \psi (2)_{i}^{n}\,|\,%
\mathcal{F}_{\frac{i-1}{n}}\right] \,\,/\sqrt{n}
\end{equation*}%
is a linear combination of the following three terms
\begin{equation*}
\mathrm{E}\left[ g^{\prime }(\beta _{i}^{n})\cdot \sigma _{\frac{i-1}{n}%
}^{\prime }\,\int_{(i-1)/n}^{i/n}\left( W_{u}-W_{\frac{i-1}{n}}\right) \,%
\mathrm{d}W_{u}\,|\,\mathcal{F}_{\frac{i-1}{n}}\right] \,,
\end{equation*}%
\begin{equation*}
\mathrm{E}\left[ g^{\prime }(\beta _{i}^{n})\cdot \triangle
_{i}^{n}M(n,i)\cdot \triangle _{i}^{n}W\,|\,\mathcal{F}_{\frac{i-1}{n}}%
\right] \,
\end{equation*}%
and%
\begin{equation*}
\mathrm{E}[\,g^{\prime }(\beta _{i}^{n})\cdot \int_{(i-1)/n}^{i/n}W_{u}\,%
\mathrm{d}M(n,i)_{u}\,|\,\mathcal{F}_{\frac{i-1}{n}}\,].
\end{equation*}%
But these three terms are all equal to $0$ as seen by the following
arguments.\vspace{1mm}

The conditional distribution of
\begin{equation*}
\left( W_{t}-W_{\frac{i-1}{n}}\right) _{t\geq \frac{i-1}{n}}|\mathcal{F}_{%
\frac{i-1}{n}}
\end{equation*}%
is clearly not affected by a change of sign. Thus since $g$ being assumed
even and $g^{\prime }$ therefore odd we have
\begin{equation*}
\mathrm{E}\left[ \,g^{\prime }(\beta _{i}^{n})\,\int_{(i-1)/n}^{i/n}\left(
W_{u}-W_{\frac{i-1}{n}}\right) \,\mathrm{d}W_{u}\,|\,\mathcal{F}_{\frac{i-1}{%
n}}\,\right] =0
\end{equation*}%
implying the vanishing of the first term. \vspace{1mm}

Secondly, by assumption, $\left( W_{t}-W_{\frac{i-1}{n}}\right) _{t\geq
\frac{i-1}{n}}$ and $(M(n,i)_{t})_{t\geq \frac{i-1}{n}}$ are independent
given $\mathcal{F}_{\frac{i-1}{n}}$. Therefore, denoting by $\mathcal{F}%
_{i,n}^{\,0}$ the $\sigma $-field generated by
\begin{equation*}
\left( W_{t}-W_{\frac{i-1}{n}}\right) _{\frac{i-1}{n}\leq t\leq i/n}\ \ \
\text{and}\ \ \ \mathcal{F}_{\frac{i-1}{n}},
\end{equation*}%
the martingale property of $(M(n,i)_{t})$ ensures that
\begin{equation*}
\mathrm{E}[\,g^{\prime }(\beta _{i}^{n})\cdot \triangle _{i}^{n}M(n,i)\cdot
\triangle _{i}^{n}W\,|\,\mathcal{F}_{i,n}^{\,0}\,]=0\
\end{equation*}%
and%
\begin{equation*}
\mathrm{E}[\left[ g^{\prime }(\beta _{i}^{n})\cdot
\int_{(i-1)/n}^{i/n}W_{u}\,\mathrm{d}M(n,i)_{u}\,|\,\mathcal{F}_{i,n}^{\,0}%
\right] \,=0.
\end{equation*}%
Using this the vanishing of
\begin{equation*}
\mathrm{E}\left[ \,g^{\prime }(\beta _{i}^{n})\cdot \triangle
_{i}^{n}M(n,i)\cdot \triangle _{i}^{n}W\,|\,\mathcal{F}_{\frac{i-1}{n}}%
\right]
\end{equation*}%
and%
\begin{equation*}
\mathrm{E}\left[ \,g^{\prime }(\beta _{i}^{n})\cdot
\int_{(i-1)/n}^{i/n}W_{u}\,\mathrm{d}M(n,i)_{u}\,|\,\mathcal{F}_{\frac{i-1}{n%
}}\right] \,\vspace{1mm}
\end{equation*}%
is easily obtained by successive conditioning.\vspace{1mm}

The proof of (\ref{13a}) is hereby completed.

\noindent $\square $

\bibliographystyle{chicago}
\bibliography{neil}

\begin{thebibliography}{}

\bibitem[\protect\citeauthoryear{Andersen and Bollerslev}{Andersen and
  Bollerslev}{1997}]{AndersenBollerslev(97jef)}
Andersen, T.~G. and T.~Bollerslev (1997).
\newblock Intraday periodicity and volatility persistence in financial markets.
\newblock {\em Journal of Empirical Finance\/}~{\em 4}, 115--158.

\bibitem[\protect\citeauthoryear{Andersen and Bollerslev}{Andersen and
  Bollerslev}{1998}]{AndersenBollerslev(98)}
Andersen, T.~G. and T.~Bollerslev (1998).
\newblock Deutsche mark-dollar volatility: intraday activity patterns,
  macroeconomic announcements, and longer run dependencies.
\newblock {\em Journal of Finance\/}~{\em 53}, 219--265.

\bibitem[\protect\citeauthoryear{Andersen, Bollerslev, and Diebold}{Andersen
  et~al.}{2003}]{AndersenBollerslevDiebold(03bipower)}
Andersen, T.~G., T.~Bollerslev, and F.~X. Diebold (2003).
\newblock Some like it smooth, and some like it rough: untangling continuous
  and jump components in measuring, modeling and forecasting asset return
  volatility.
\newblock Unpublished paper: Economics Dept, Duke University.

\bibitem[\protect\citeauthoryear{Andersen, Bollerslev, and Diebold}{Andersen
  et~al.}{2005}]{AndersenBollerslevDiebold(05)}
Andersen, T.~G., T.~Bollerslev, and F.~X. Diebold (2005).
\newblock Parametric and nonparametric measurement of volatility.
\newblock In Y.~A{\"{\i}}t-Sahalia and L.~P. Hansen (Eds.), {\em Handbook of
  Financial Econometrics}. Amsterdam: North Holland.
\newblock Forthcoming.

\bibitem[\protect\citeauthoryear{Andersen, Bollerslev, Diebold, and
  Labys}{Andersen et~al.}{2001}]{AndersenBollerslevDieboldLabys(01)}
Andersen, T.~G., T.~Bollerslev, F.~X. Diebold, and P.~Labys (2001).
\newblock The distribution of exchange rate volatility.
\newblock {\em Journal of the American Statistical Association\/}~{\em 96},
  42--55.
\newblock Correction published in 2003, volume 98, page 501.

\bibitem[\protect\citeauthoryear{Andersen, Bollerslev, Diebold, and
  Labys}{Andersen et~al.}{2003}]{AndersenBollerslevDieboldLabys(03model)}
Andersen, T.~G., T.~Bollerslev, F.~X. Diebold, and P.~Labys (2003).
\newblock Modeling and forecasting realized volatility.
\newblock {\em Econometrica\/}~{\em 71}, 579--625.

\bibitem[\protect\citeauthoryear{Bandi and Russell}{Bandi and
  Russell}{2003}]{BandiRussell(03)}
Bandi, F.~M. and J.~R. Russell (2003).
\newblock Microstructure noise, realized volatility, and optimal sampling.
\newblock Unpublished paper, Graduate School of Business, University of
  Chicago.

\bibitem[\protect\citeauthoryear{Barndorff-Nielsen, Graversen, Jacod,
  Podolskij, and Shephard}{Barndorff-Nielsen
  et~al.}{2004}]{BarndorffNielsenGraversenJacodPodolskyShephard(04shiryaev)}
Barndorff-Nielsen, O.~E., S.~E. Graversen, J.~Jacod, M.~Podolskij, and
  N.~Shephard (2004).
\newblock A central limit theorem for realised power and bipower variations of
  continuous semimartingales.
\newblock Economics working paper 2004-W29, Nuffield College, Oxford.

\bibitem[\protect\citeauthoryear{Barndorff-Nielsen, Hansen, Lunde, and
  Shephard}{Barndorff-Nielsen
  et~al.}{2004}]{BarndorffNielsenHansenLundeShephard(04)}
Barndorff-Nielsen, O.~E., P.~R. Hansen, A.~Lunde, and N.~Shephard (2004).
\newblock Regular and modified kernel-based estimators of integrated variance:
  the case with independent noise.
\newblock Unpublished paper: Nuffield College, Oxford.

\bibitem[\protect\citeauthoryear{Barndorff-Nielsen and
  Shephard}{Barndorff-Nielsen and
  Shephard}{2002}]{BarndorffNielsenShephard(02realised)}
Barndorff-Nielsen, O.~E. and N.~Shephard (2002).
\newblock Econometric analysis of realised volatility and its use in estimating
  stochastic volatility models.
\newblock {\em Journal of the Royal Statistical Society, Series B\/}~{\em 64},
  253--280.

\bibitem[\protect\citeauthoryear{Barndorff-Nielsen and
  Shephard}{Barndorff-Nielsen and
  Shephard}{2003}]{BarndorffNielsenShephard(03bernoulli)}
Barndorff-Nielsen, O.~E. and N.~Shephard (2003).
\newblock Realised power variation and stochastic volatility.
\newblock {\em Bernoulli\/}~{\em 9}, 243--265.
\newblock Correction published in pages 1109--1111.

\bibitem[\protect\citeauthoryear{Barndorff-Nielsen and
  Shephard}{Barndorff-Nielsen and
  Shephard}{2004a}]{BarndorffNielsenShephard(04multi)}
Barndorff-Nielsen, O.~E. and N.~Shephard (2004a).
\newblock Econometric analysis of realised covariation: high frequency
  covariance, regression and correlation in financial economics.
\newblock {\em Econometrica\/}~{\em 72}, 885--925.

\bibitem[\protect\citeauthoryear{Barndorff-Nielsen and
  Shephard}{Barndorff-Nielsen and
  Shephard}{2004b}]{BarndorffNielsenShephard(04jfe)}
Barndorff-Nielsen, O.~E. and N.~Shephard (2004b).
\newblock Power and bipower variation with stochastic volatility and jumps
  (with discussion).
\newblock {\em Journal of Financial Econometrics\/}~{\em 2}, 1--48.

\bibitem[\protect\citeauthoryear{Barndorff-Nielsen and
  Shephard}{Barndorff-Nielsen and
  Shephard}{2005a}]{BarndorffNielsenShephard(03test)}
Barndorff-Nielsen, O.~E. and N.~Shephard (2005a).
\newblock Econometrics of testing for jumps in financial economics using
  bipower variation.
\newblock {\em Journal of Financial Econometrics\/}.
\newblock Forthcoming.

\bibitem[\protect\citeauthoryear{Barndorff-Nielsen and
  Shephard}{Barndorff-Nielsen and
  Shephard}{2005b}]{BarndorffNielsenShephard(05tom)}
Barndorff-Nielsen, O.~E. and N.~Shephard (2005b).
\newblock How accurate is the asymptotic approximation to the distribution of
  realised volatility?
\newblock In D.~W.~K. Andrews, J.~Powell, P.~A. Ruud, and J.~H. Stock (Eds.),
  {\em Identification and Inference for Econometric Models. {A} Festschrift in
  Honour of {T.J. Rothenberg}}. Cambridge: Cambridge University Press.
\newblock Forthcoming.

\bibitem[\protect\citeauthoryear{Barndorff-Nielsen, Shephard, and
  Winkel}{Barndorff-Nielsen et~al.}{2004}]{BarndorffNielsenShephardWinkel(04)}
Barndorff-Nielsen, O.~E., N.~Shephard, and M.~Winkel (2004).
\newblock Limit theorems for multipower variation in the presence of jumps in
  financial econometrics.
\newblock Unpublished paper: Nuffield College, Oxford.

\bibitem[\protect\citeauthoryear{Calvet and Fisher}{Calvet and
  Fisher}{2002}]{CalvetFisher(02)}
Calvet, L. and A.~Fisher (2002).
\newblock Multifractality in asset returns: theory and evidence.
\newblock {\em Review of Economics and Statistics\/}~{\em 84}, 381--406.

\bibitem[\protect\citeauthoryear{Corradi and Distaso}{Corradi and
  Distaso}{2004}]{CorradiDistaso(04)}
Corradi, V. and W.~Distaso (2004).
\newblock Specification tests for daily integrated volatility, in the presence
  of possible jumps.
\newblock Unpublished paper: Queen Mary College, London.

\bibitem[\protect\citeauthoryear{Delattre and Jacod}{Delattre and
  Jacod}{1997}]{DelattreJacod(97)}
Delattre, S. and J.~Jacod (1997).
\newblock A central limit theorem for normalized functions of the increments of
  a diffusion process in the presence of round off errors.
\newblock {\em Bernoulli\/}~{\em 3}, 1--28.

\bibitem[\protect\citeauthoryear{Doob}{Doob}{1953}]{Doob(53)}
Doob, J.~L. (1953).
\newblock {\em Stochastic Processes}.
\newblock New York: John Wiley and Sons.

\bibitem[\protect\citeauthoryear{Forsberg and Ghysels}{Forsberg and
  Ghysels}{2004}]{ForsbergGhysels(04)}
Forsberg, L. and E.~Ghysels (2004).
\newblock Why do absolute returns predict volatility so well.
\newblock Unpublished paper: Economics Department, UNC, Chapel Hill.

\bibitem[\protect\citeauthoryear{Ghysels, Harvey, and Renault}{Ghysels
  et~al.}{1996}]{GhyselsHarveyRenault(96)}
Ghysels, E., A.~C. Harvey, and E.~Renault (1996).
\newblock Stochastic volatility.
\newblock In C.~R. Rao and G.~S. Maddala (Eds.), {\em Statistical Methods in
  Finance}, pp.\  119--191. Amsterdam: North-Holland.

\bibitem[\protect\citeauthoryear{Ghysels, Santa-Clara, and Valkanov}{Ghysels
  et~al.}{2004}]{GhyselsSantaClaraValkoanov(04)}
Ghysels, E., P.~Santa-Clara, and R.~Valkanov (2004).
\newblock Predicting volatility: getting the most out of return data sampled at
  different frequencies.
\newblock Unpublished paper: Department of Economics, University of North
  Carolina.

\bibitem[\protect\citeauthoryear{Gloter and Jacod}{Gloter and
  Jacod}{2001a}]{GloterJacod(01a)}
Gloter, A. and J.~Jacod (2001a).
\newblock Diffusions with measurement errors. \uppercase{I} --- local
  asymptotic normality.
\newblock {\em ESAIM: Probability and Statistics\/}~{\em 5}, 225--242.

\bibitem[\protect\citeauthoryear{Gloter and Jacod}{Gloter and
  Jacod}{2001b}]{GloterJacod(01b)}
Gloter, A. and J.~Jacod (2001b).
\newblock Diffusions with measurement errors. \uppercase{II} --- measurement
  errors.
\newblock {\em ESAIM: Probability and Statistics\/}~{\em 5}, 243--260.

\bibitem[\protect\citeauthoryear{Goncalves and Meddahi}{Goncalves and
  Meddahi}{2004}]{GoncalvesMeddahi(04)}
Goncalves, S. and N.~Meddahi (2004).
\newblock Bootstrapping realized volatility.
\newblock Unpublished paper, CIRANO, Montreal.

\bibitem[\protect\citeauthoryear{Hansen and Lunde}{Hansen and
  Lunde}{2003}]{HansenLunde(03)}
Hansen, P.~R. and A.~Lunde (2003).
\newblock An optimal and unbiased measure of realized variance based on
  intermittent high-frequency data.
\newblock Unpublished paper, Department of Economics, Stanford University.

\bibitem[\protect\citeauthoryear{Huang and Tauchen}{Huang and
  Tauchen}{2003}]{HuangTauchen(03)}
Huang, X. and G.~Tauchen (2003).
\newblock The relative contribution of jumps to total price variation.
\newblock Unpublished paper: Department of Economics, Duke University.

\bibitem[\protect\citeauthoryear{Jacod}{Jacod}{1994}]{Jacod(94)}
Jacod, J. (1994).
\newblock Limit of random measures associated with the increments of a
  \uppercase{B}rownian semimartingale.
\newblock Preprint number 120, Laboratoire de Probabiliti\'{e}s, Universit\'{e}
  Pierre et Marie Curie, Paris.

\bibitem[\protect\citeauthoryear{Jacod, Lejay, and Talay}{Jacod
  et~al.}{2005}]{JacodLejayTalay(05)}
Jacod, J., A.~Lejay, and D.~Talay (2005).
\newblock Testing the multiplicity of a diffusion.
\newblock In preparation.

\bibitem[\protect\citeauthoryear{Jacod and Protter}{Jacod and
  Protter}{1998}]{JacodProtter(98)}
Jacod, J. and P.~Protter (1998).
\newblock Asymptotic error distributions for the \uppercase{E}uler method for
  stochastic differential equations.
\newblock {\em Annals of Probability\/}~{\em 26}, 267--307.

\bibitem[\protect\citeauthoryear{Jacod and Shiryaev}{Jacod and
  Shiryaev}{2003}]{JacodShiryaev(03)}
Jacod, J. and A.~N. Shiryaev (2003).
\newblock {\em Limit Theorems for Stochastic Processes\/} (2 ed.).
\newblock Springer-Verlag: Berlin.

\bibitem[\protect\citeauthoryear{Karatzas and Shreve}{Karatzas and
  Shreve}{1991}]{KaratzasShreve(91)}
Karatzas, I. and S.~E. Shreve (1991).
\newblock {\em Brownian Motion and Stochastic Calculus\/} (2 ed.), Volume 113
  of {\em Graduate Texts in Mathematics}.
\newblock Berlin: Springer--Verlag.

\bibitem[\protect\citeauthoryear{Karatzas and Shreve}{Karatzas and
  Shreve}{1998}]{KaratzasShreve(98)}
Karatzas, I. and S.~E. Shreve (1998).
\newblock {\em Methods of Mathematical Finance}.
\newblock New York: Springer--Verlag.

\bibitem[\protect\citeauthoryear{Maheswaran and Sims}{Maheswaran and
  Sims}{1993}]{MaheswaranSims(93)}
Maheswaran, S. and C.~A. Sims (1993).
\newblock Empirical implications of arbitrage-free asset markets.
\newblock In P.~C.~B. Phillips (Ed.), {\em Models, Methods and Applications of
  Econometrics}, pp.\  301--316. Basil Blackwell.

\bibitem[\protect\citeauthoryear{Munroe}{Munroe}{1953}]{Munroe(53)}
Munroe, M.~E. (1953).
\newblock {\em Introduction to Measure and Integration}.
\newblock Cambridge, MA: Addison-Wesley Publishing Company, Inc.

\bibitem[\protect\citeauthoryear{Mykland and Zhang}{Mykland and
  Zhang}{2005}]{MyklandZhang(05)}
Mykland, P. and L.~Zhang (2005).
\newblock \uppercase{ANOVA} for diffusions.
\newblock {\em Annals of Statistics\/}~{\em 33}.
\newblock Forthcoming.

\bibitem[\protect\citeauthoryear{Nielsen and Frederiksen}{Nielsen and
  Frederiksen}{2005}]{NielsenFrederiksen(05)}
Nielsen, M.~O. and P.~H. Frederiksen (2005).
\newblock Finite sample accuracy of integrated volatility estimators.
\newblock Unpublished paper, Department of Economics, Cornell University.

\bibitem[\protect\citeauthoryear{Parkinson}{Parkinson}{1980}]{Parkinson(80)}
Parkinson, M. (1980).
\newblock The extreme value method for estimating the variance of the rate of
  return.
\newblock {\em Journal of Business\/}~{\em 53}, 61--66.

\bibitem[\protect\citeauthoryear{Revuz and Yor}{Revuz and
  Yor}{1999}]{RevuzYor(99)}
Revuz, D. and M.~Yor (1999).
\newblock {\em Continuous Martingales and Brownian motion\/} (3 ed.).
\newblock Heidelberg: Springer-Verlag.

\bibitem[\protect\citeauthoryear{Schwert}{Schwert}{1990}]{Schwert(90JB)}
Schwert, G.~W. (1990).
\newblock Indexes of \uppercase{U.S.} stock prices from 1802 to 1987.
\newblock {\em Journal of Business\/}~{\em 63}, 399--426.

\bibitem[\protect\citeauthoryear{Shephard}{Shephard}{2005}]{Shephard(05)}
Shephard, N. (2005).
\newblock {\em Stochastic Volatility: Selected Readings}.
\newblock Oxford: Oxford University Press.
\newblock Forthcoming.

\bibitem[\protect\citeauthoryear{Shiryaev}{Shiryaev}{1999}]{Shiryaev(99)}
Shiryaev, A.~N. (1999).
\newblock {\em Essentials of Stochastic Finance: Facts, Models and Theory}.
\newblock Singapore: World Scientific.

\bibitem[\protect\citeauthoryear{Woerner}{Woerner}{2004}]{Woerner(04power)}
Woerner, J. (2004).
\newblock Power and multipower variation: inference for high frequency data.
\newblock Unpublished paper.

\bibitem[\protect\citeauthoryear{Zhang}{Zhang}{2004}]{Zhang(04)}
Zhang, L. (2004).
\newblock Efficient estimation of stochastic volatility using noisy
  observations: a multi-scale approach.
\newblock Unpublished paper: Department of Statistics, Carnegie Mellon
  University.

\bibitem[\protect\citeauthoryear{Zhang, Mykland, and A{\"{\i}}t-Sahalia}{Zhang
  et~al.}{2005}]{ZhangMyklandAitSahalia(03)}
Zhang, L., P.~Mykland, and Y.~A{\"{\i}}t-Sahalia (2005).
\newblock A tale of two time scales: determining integrated volatility with
  noisy high-frequency data.
\newblock {\em Journal of the American Statistical Association\/}.
\newblock Forthcoming.

\end{thebibliography}

\end{document}